\documentclass[onefignum,onetabnum]{siamonline250211}

\usepackage{amsmath,amssymb,mathtools}
\usepackage{graphicx}
\usepackage{placeins}
\usepackage{adjustbox}
\usepackage{booktabs,array}
\usepackage{enumitem}
\usepackage{xcolor}

\newsiamthm{conjecture}{Conjecture}
\newsiamremark{remark}{Remark}
\newsiamremark{empirical}{Empirical statement}
\newsiamremark{example}{Example}                            %
\crefname{conjecture}{Conjecture}{Conjectures}
\crefname{empirical}{Empirical statement}{Empirical statements}
\crefname{example}{Example}{Examples}

\newcommand{\RR}{\mathbb{R}}
\newcommand{\TT}{\mathbb{T}}
\providecommand{\rank}{\operatorname{rank}}
\newcommand{\rch}{\operatorname{reach}}
\providecommand{\tr}{\operatorname{tr}}
\newcommand{\II}{\mathrm{I\!I}}
\newcommand{\normal}{\mathbf{n}}

\newcommand{\folord}{\operatorname{ord}}                  %

\headers{Coordinate dependence of numerical rank}{H. Tran and D. P. Truong}

\title{The coordinate dependence of numerical rank: charts adapted to level sets%
\thanks{The authors contributed equally.}}

\author{Hung Tran\thanks{Department of Mathematics and Statistics, Texas Tech University, Lubbock, TX (\email{hung.tran@ttu.edu}).}
\and Duc P. Truong\thanks{Theoretical Division, Los Alamos National Laboratory, Los Alamos, NM (\email{dptruong@lanl.gov}).}}

\ifpdf
\hypersetup{pdftitle={The coordinate dependence of numerical rank: charts adapted to level sets}}
\fi

\begin{document}
\maketitle

\begin{abstract}
Sampling a smooth field on a grid can give an array of high numerical rank even when the field is simple. A Gaussian ridge
rotated off the grid axes is the standard example: in coordinates which follow the ridge, the same field has rank one. The
rank therefore depends on the chart, that is, on the coordinates in which the field is sampled. We ask which properties of a
chart control this rank.
For a field which depends on position only through a function $\phi$, a chart whose coordinate lines follow the level sets
of $\phi$ gives rank one. Building on this, we show four things. First, mesh quality does not control rank: a chart whose
cells have perfect size and shape can still give infinite rank. Second, when a chart follows the level sets of an estimate
of $\phi$, the extra rank depends on a single quantity, the defect, which measures how much the true field varies along each
coordinate line. Its size matters only logarithmically, its complexity through a count which cannot be improved in general,
and when the defect is simple across the level sets the cost stays small. Third, whether the defect is simple can be
checked by differentiating the true level function (or, for a monotone profile, the field itself) along the flow which
builds the chart. Fourth, rank and invertibility fail differently. When the curve behind an offset chart is perturbed, the rank grows
slowly, but the chart stops being invertible abruptly, once the offset exceeds the distance at which neighboring normals
of the curve first cross. A chart built from the level sets of a field estimate stays invertible when the gradient of the
estimate does not vanish on the band and the coordinate lines along each level set do not cross themselves.
Pre-registered experiments in two dimensions show how far the bounds are from the observed ranks, and suggest when each kind
of chart is preferable; the tolerance behind the last comparison, and a correction to how it was scored, were chosen after
the experiments.
\end{abstract}



\section{Introduction}\label{sec:intro}\label{note:sec:intro}\label{note:sec:lit}

Sampling a smooth field $U$ on a Cartesian grid can produce an array of high numerical rank, even when $U$ is simple. The standard example is a Gaussian ridge rotated off the coordinate axes. On the other hand, in a chart adapted to
the level sets of the ridge, that is, a map from a computational domain onto the physical domain whose coordinate lines
follow those level sets, the same field gives an array of rank one. Therefore, the numerical rank of a sampled field is a
property of the field together with its chart, not of the field alone. This work grew out of tensor-train solvers for
partial differential equations~\cite{DanisWENO2025,TranIGA2026}, in which the ranks of solutions and of
geometric coefficients depend strongly on the grid or mapping on which they are represented.

\begin{figure}[tbhp]
\centering
\includegraphics{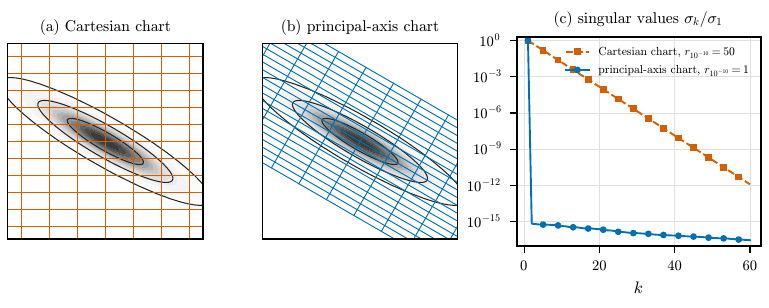}
\caption{Rank depends on the chart. The rotated Gaussian $U=\exp(-x^{\mathsf T}Ax)$, $A=Q^{\mathsf T}\mathrm{diag}(1,\kappa_\Sigma)Q$,
$\kappa_\Sigma=25$, rotated by $30^\circ$, is sampled on $201\times201$ nodes of (a) the Cartesian chart and (b) the chart along
the principal axes.
(c) Normalized singular values of the two sampled matrices, with their relative Frobenius $\varepsilon$-ranks at
$\varepsilon=10^{-10}$.}\label{fig:intro}
\end{figure}

A rotated Gaussian makes the dependence on the chart quantitative. For the variance ratio $\kappa_\Sigma=25$ and a rotation by
$30^\circ$, the Cartesian $\varepsilon$-rank at $\varepsilon=10^{-10}$ is $50$. In the chart along the principal axes, the
Gaussian is a product of two functions of one variable each, so its $\varepsilon$-rank is $1$ (Figure~\ref{fig:intro}). Both charts are isometries with condition number one, so the difference in rank is not a matter of
the chart's metric. In this paper, we ask \emph{which geometric properties of a chart $X$ control the rank of $U\circ X$, and
quantitatively how}. %
Our answer is in terms of the chart's foliation, that is, of how its coordinate surfaces sit relative to the level sets of
the field. If they lie in level sets, the rank is one; this is sufficient, not necessary (Section~\ref{sec:condnumber}). If
they do not, we bound the extra rank by how far they are from doing so. The metric criteria of mesh adaptation do not
determine the rank.

This paper makes four contributions, one in each of Sections~\ref{sec:metric}--\ref{sec:invert}.
\begin{enumerate}[label=(\arabic*),leftmargin=*]
\item \emph{Mesh quality does not control rank} (Section~\ref{sec:metric}). Mesh adaptation measures the quality of a
chart by the shape and size of its cells, through two standard criteria, alignment and
equidistribution~\cite{HuangRussell2011}. We show that a chart of perfect quality, which meets both criteria exactly, can
still give a field of infinite rank (Theorem~\ref{thm:counterexample}).
\item \emph{The rank cost of misalignment} (Section~\ref{sec:misalign}). In practice, a chart follows the level sets of an
estimate of the field, not of the field itself, so it is misaligned. We show that the extra rank depends on one quantity,
the \emph{defect}: how much the true field varies along each coordinate line of the chart. Its size matters only
logarithmically (Theorem~\ref{thm:rankcost}). Its complexity matters more, through a count which cannot be improved in
general (Proposition~\ref{prop:binomialattained}). When the defect is simple across the level sets, for instance
polynomial or analytic in the level coordinate, the cost stays small (Theorem~\ref{thm:productclosed} and
Corollary~\ref{cor:levelanalytic}).
\item \emph{The geometric source of the cost} (Section~\ref{sec:foliation}). Whether the defect is simple can be checked
by differentiating the true level function along the curves which build the chart (Proposition~\ref{prop:liederivs}), or
the field itself when the profile is monotone (Remark~\ref{rem:folorderinvariance}). For example, the chart obtained by
offsetting a perturbed straight line has a defect with just two terms (Corollary~\ref{cor:orderstraight}). For curved
charts, the cost is bounded in terms of the width of the band compared with the distance to the nearest point where the
defect stops being analytic (Corollary~\ref{cor:georank}).
\item \emph{Rank and invertibility fail differently} (Section~\ref{sec:invert}). To reconstruct the field, the chart must
also be invertible. For a chart built by offsetting a perturbed straight line, the perturbation raises the rank only
slowly, but the chart stops being injective once the offset exceeds the reach of the perturbed curve, which is at most its
smallest radius of curvature (Theorem~\ref{thm:rankinjdichotomy}). A chart built from the level sets of a field estimate
has no such limit: it is injective whenever the estimate has no critical points on the band and its leaf charts are
injective (Theorem~\ref{thm:leafinjectivity}). Finally, experiments suggest when
each kind of chart is preferable (Empirical statement~\ref{emp:chartchoice}), with the caveat that the tolerance behind
this, and a correction of a coverage test, were chosen after the experiments (Section~\ref{sec:numerics}).
\end{enumerate}

The rank bound of Theorem~\ref{thm:rankcost} follows the mechanism of M.~Udell and A.~Townsend~\cite{UdellTownsend2019},
which expands the generating function and counts the resulting separable terms. We use a Chebyshev expansion~\cite{Trefethen2019} and the
relative Frobenius norm. The generic rank of Hadamard powers is due to~\cite{DammDietrich2023}, and the count is attained
for zero shifts, when all Taylor coefficients of $G$ up to degree $m$ are nonzero, by~\cite[Thm.~3.2]{MaddenThrampoulidis2024}. %
What is new is the reduction of the rank cost of a misaligned chart to its defect, the attainment of the count for generic
level shifts, the bounds for simple defects, and their characterization by derivatives along the flow which builds the
chart. The other classical tools are credited where they are used.

\emph{Related work.} Several works choose coordinates to lower tensor rank, by optimizing over linear changes of
coordinates~\cite{DektorVenturi2023,DektorVenturi2025}, or by learning invertible maps which are composed with tensor
trains~\cite{Eigel2025}. A related approach writes a field as a fixed profile of a low-rank function~\cite{Krah2023}.
These works construct a good map. We ask instead which properties of a given map control the rank. In model reduction and
shock tracking, maps are chosen to make a family of solutions compressible~\cite{Taddei2020,Iollo2026}, or to align a mesh
with a discontinuity~\cite{ZahrPersson2018,MirhoseiniZahr2023} (see the survey~\cite{Hesthaven2026}). These are
approximation goals for a family or for accuracy, not the rank of one field. Bounds on the rank of matrices sampled from
smooth functions~\cite{UdellTownsend2019,BigDataLowRank2025,TownsendTrefethen2015,Griebel2023,BeckermannTownsend2017} take
the sampling as fixed, and do not ask how a change of coordinates affects it. Mesh
adaptation~\cite{HuangRussell2011,Loseille2011,HuangKamenski2018}, optimal-transport meshes~\cite{BuddWilliams2009}, and
surface parametrization~\cite{Ray2008,Kalberer2007,Bommes2009,Vaxman2016} build charts with well-shaped cells, but not
with rank as the goal, and Section~\ref{sec:metric} shows why that is not enough. Finally, the curves which build our
charts also appear in active manifolds~\cite{Bridges2019} and in level-set learning~\cite{ZhangZhangLin2019,Teng2023}.

\emph{Scope and outline.} The bounds concern the (fiber)\,$\times$\,(level) unfolding of the sampled field, which is the
whole question for $d=2$ and one tensor-train unfolding for $d\ge3$. The experiments are two-dimensional. The three experiments which serve as evidence were
pre-registered (Section~\ref{sec:numerics}). Section~\ref{sec:setting} recalls the exact theory, and
Sections~\ref{sec:metric}--\ref{sec:invert} contain the contributions (1)--(4). Section~\ref{sec:numerics} reports the
experiments, Section~\ref{sec:discussion} discusses the cost of the chart, the limits of a single chart, open problems, and
limitations, and Appendix~\ref{app:proofs} contains the longer proofs.

\section{The exact theory}\label{sec:setting}\label{note:sec:setting}\label{note:sec:exact}

Here, we fix the representation, the measure of rank, and the exact theory which the rest of the paper perturbs.
The results here are classical or elementary. %

Let $\Omega\subset\RR^d$ be open and let $U:\Omega\to\RR$. A chart is a map
$X:\Omega_c\to\Omega$ defined on a computational domain. The geometry-aware representation is the pair
\[
\mathcal{R}=\bigl(X,\ \widetilde U\bigr),\qquad \widetilde U=U\circ X ,
\]
and the reconstruction is the \emph{inverse-coordinate pullback}
\begin{equation}\label{note:eq:pullback}
\widehat u(x)=\widetilde U_r\bigl(X^{-1}(x)\bigr),
\end{equation}
where $\widetilde U_r$ is a rank-$r$ approximation of $\widetilde U$.

Throughout, we employ~\eqref{note:eq:pullback}, rather than the scattered interpolation of the values of
$\widetilde U_r$ at the physical points $X(\xi)$. At a
prescribed reconstruction error, the cost of the representation is $C_X+C_{\widetilde U}+C_{X^{-1}}$, that is, the storage
of the chart, of the factors of $\widetilde U_r$, and of whatever evaluating $X^{-1}$ requires. Thus, storing the chart
counts, and Section~\ref{sec:discussion} returns to this cost.

\emph{Notation.} We write the computational coordinates as $\xi=(\xi',\xi_d)$ with $\xi'=(\xi_1,\dots,\xi_{d-1})$. In a
chart, $\xi'$ runs along the level sets of a function and $\xi_d$ across them: we call $\xi'$
the fiber coordinate and $\xi_d$ the level coordinate. We call the level sets of a function $\phi$ its \emph{leaves}, and
the family of all of them the \emph{foliation} of $\phi$. Likewise, the coordinate surfaces $\{\xi_d=t\}$ of a chart are
its leaves, and together they form the \emph{foliation of the chart}. Throughout, a sampled field is a matrix whose rows are indexed by
fiber points and whose columns are indexed by levels (fiber $\times$ level).

\begin{definition}\label{def:epsrank}
Throughout, the relative Frobenius $\varepsilon$-rank of a matrix $\mathbf A$ is
\[
r_\varepsilon(\mathbf A)=\min\{r:\|\mathbf A-\mathbf A_r\|_F\le\varepsilon\|\mathbf A\|_F\},
\]
where $\mathbf A_r$ is a best rank-$r$ approximation of $\mathbf A$. %
For a kernel $f$, $r_\varepsilon(f)$ is defined in the same way, with the Hilbert--Schmidt norm.
\end{definition}

\subsection{The exact theory on a regular band}

If $U(x)=G(\phi(x))$, then $U$ is constant on each level set of $\phi$. In a chart whose coordinate lines follow the
level sets, $U$ therefore depends on the level coordinate alone, and its sampled matrix has rank one. Such a chart exists
on any compact band where $\nabla\phi\ne0$, and a flow builds it.

Let $\phi\in C^r(\Omega)$, $r\ge2$, on an open set $\Omega\subset\RR^d$, and let the band $\mathcal B=\phi^{-1}([a,b])$ be
nonempty and compact, with $\nabla\phi\ne0$ on $\mathcal B$. The flow $\Phi_t$ of $V=\nabla\phi/|\nabla\phi|^2$ raises the
level by $t$: $\phi(\Phi_t(x))=\phi(x)+t$.

\begin{theorem}[Regular-band structure]\label{note:thm:band}\label{thm:bandproduct}
Let $c_\ast\in(a,b)$ and $F:=\phi^{-1}(c_\ast)$. Then the map
\[
X:F\times[a,b]\to\mathcal B,\qquad X(y,c):=\Phi_{c-c_\ast}(y)
\]
is a $C^{r-1}$ diffeomorphism with $\phi(X(y,c))=c$. Consequently, $\mathcal B\cong F\times[a,b]$. %
\end{theorem}

\begin{figure}[tbhp]
\centering
\includegraphics{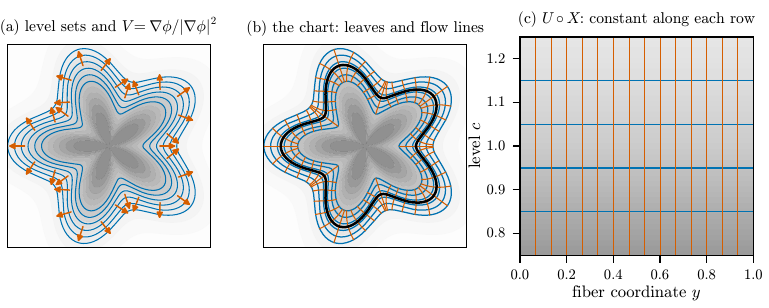}
\caption{The construction of Theorem~\ref{thm:bandproduct} on the wavy radial field $U=\exp(-\phi^2)$,
$\phi=|x|+a_w\cos(K\vartheta)$, $a_w=0.25$, $K=5$, on the band of levels $[0.75,1.25]$. (a) Level sets and the field
$V=\nabla\phi/|\nabla\phi|^2$ on the leaf $\phi=1$. (b) The chart. The flow of $V$ carries the leaf $\phi=1$ (black) onto
the other leaves, and the leaves (blue) and flow lines (red) form the coordinate mesh. (c) $U\circ X$ on the computational
rectangle (fiber coordinate $\times$ level). It is constant in the fiber coordinate, that is, on each leaf, hence of rank
one (Proposition~\ref{prop:rankone}).}\label{fig:construction}
\end{figure}

That is, the flow carries one level set $F$ onto all the others, so a point of the band is described by its starting
point $y$ on $F$ and its level $c$. This is the proof of Milnor's regular-interval theorem~\cite[Thm.~3.1]{Milnor1963};
see also Ehresmann's fibration theorem~\cite{Ehresmann1951}. Section~\ref{sec:foliation} builds charts by the same flow
from an estimate of $\phi$. Figure~\ref{fig:construction} shows the construction on the wavy radial field of
Section~\ref{sec:numerics}.

\begin{proposition}[Exact rank collapse]\label{note:prop:collapse}\label{prop:rankone}
Suppose that $U=G\circ\phi$. Then $(U\circ X)(y,c)=G(c)$. Consequently, for samples $y_i\in F$, $c_j\in[a,b]$ with $\mathbf g_j=G(c_j)\ne0$ for some $j$,
\[\mathbf A_{ij}=U(X(y_i,c_j))=\mathbf{1}\,\mathbf g^{\mathsf T},\]
so $\rank\mathbf A=1$. Furthermore, let ${Y}:D_1\times\cdots\times D_{d-1}\to F$ be a
tensorizable fiber chart, that is, any map from a product of intervals and circles into $F$. Then the resulting tensor has
canonical polyadic (CP) rank $\rank_{\mathrm{CP}}=1$ and tensor-train (TT) rank $r_{\mathrm{TT}}=(1,\dots,1)$. %
\end{proposition}

\begin{proof}
By Theorem~\ref{note:thm:band}, $\phi(X(y,c))=c$. Thus, $U(X(y,c))=G(\phi(X(y,c)))=G(c)$, and every row of $\mathbf A$
equals $\mathbf g^{\mathsf T}$. Next, the tensor $T_{i_1\dots i_{d-1}j}=G(c_j)$ depends on the last index alone. Hence,
$T=\mathbf 1\otimes\cdots\otimes\mathbf 1\otimes\mathbf g$, which has CP rank one since $\mathbf g\ne0$. Furthermore, every TT
unfolding is the outer product of a constant vector with a nonzero vector, so it has rank one.
\end{proof}

The rank-one property does not depend on how the points are spread along each level set. We may move the points within
each leaf, differently on each leaf if we like, and every sample stays the same. We call such a choice of positions along
the leaves a \emph{gauge}.

\begin{proposition}[Gauge freedom]\label{note:prop:gauge}\label{prop:gaugefree}
Suppose that $U=G\circ\phi$, and for each level $c$ let $\pi_c:F\to F$ be any map. Then the chart
$\widehat X(y,c)=X(\pi_c(y),c)$ gives the same values as $X$: $U(\widehat X(y,c))=U(X(y,c))=G(c)$.
\end{proposition}

\begin{proof}
Since $\pi_c(y)\in F$, Theorem~\ref{note:thm:band} gives $\phi(\widehat X(y,c))=c$, so $U(\widehat X(y,c))=G(c)$.
\end{proof}

\subsection{Foliation compatibility}

The flow chart of Theorem~\ref{thm:bandproduct} is not the only chart which gives rank one.
Proposition~\ref{prop:rankone} uses only one property of it: each leaf of the chart lies inside a level set of $\phi$. The
following definition names that property.

\begin{definition}\label{note:def:fc}\label{def:foliationcompat}
A chart $X:\Omega_c\to\mathcal B$ is \emph{foliation-compatible} for $\phi$ if each of its leaves $\{\xi_d=t\}$ lies in
a level set of $\phi$, that is,
\begin{equation}\label{note:eq:fc}
\nabla\phi\cdot\partial_{\xi_i}X=0,\qquad i=1,\dots,d-1.
\end{equation}
\end{definition}

If each leaf $\{\xi_d=t\}$ of the computational domain is connected, as for a product of intervals and circles,
then~\eqref{note:eq:fc} says that $\phi\circ X$ depends on $\xi_d$ alone. If, moreover, $X$ is a diffeomorphism of
$\Omega_c=D\times I$, with $D$ a product of intervals and circles (circles are needed when the leaves are closed), onto a
band whose level sets are connected, then the leaves of the chart are exactly the level sets of $\phi$.

\begin{theorem}[Separability]\label{note:thm:fc}\label{thm:separability}
Let $X$ be foliation-compatible on a computational domain whose leaves $\{\xi_d=t\}$ are connected, let $U=G\circ\phi$,
and suppose that the sampled values of $U\circ X$ are not all zero. Then $U\circ X$ depends on the level coordinate alone,
\[
U\circ X=\bigl(G\circ\bar\Psi\bigr)(\xi_d),\qquad \bar\Psi:=\phi\circ X ,
\]
so the sampled matrix has rank one in the (fiber)\,$\times$\,(level) splitting. Under a tensorizable fiber chart
(Proposition~\ref{prop:rankone}), the sampled tensor has CP rank one and TT rank $(1,\dots,1)$.
\end{theorem}

\begin{proof}
Condition~\eqref{note:eq:fc} says that $\partial_{\xi_i}(\phi\circ X)=0$ for $i<d$. Thus, $\phi\circ X$ is constant on
each connected leaf $\{\xi_d=t\}$. That is, $\phi\circ X=\bar\Psi(\xi_d)$ and $U\circ X=G(\bar\Psi(\xi_d))$. In the
(fiber)\,$\times$\,(level) splitting, sampling gives $\mathbf A_{ij}=\mathbf g_j$ with $\mathbf g_j=G(\bar\Psi(\xi_{d,j}))$,
so $\mathbf A=\mathbf 1\mathbf g^{\mathsf T}$ has rank one, since $\mathbf g\ne0$ by assumption. Over a tensorizable fiber chart, the tensor is
$\mathbf 1\otimes\cdots\otimes\mathbf 1\otimes\mathbf g$, exactly as in the proof of Proposition~\ref{note:prop:collapse}.
\end{proof}

Foliation compatibility is sufficient for rank one, but not necessary. For example, if $G=\exp$ and $\phi\circ X=a(\xi')+b(\xi_d)$, then
$U\circ X=e^{a(\xi')}e^{b(\xi_d)}$ has rank one, while~\eqref{note:eq:fc} fails wherever $\nabla_{\xi'}a\ne0$.

Sections~\ref{sec:metric}, \ref{sec:misalign}, and~\ref{sec:invert} examine this condition from three sides: whether
the metric criteria of mesh adaptation enforce it, what it costs when it holds only approximately, as for a chart built from
an estimate of $\phi$, and what invertibility adds to it.

\section{Metric criteria do not control rank}\label{sec:metric}\label{note:sec:mm}\label{note:sec:gaussian}

Mesh adaptation designs a chart by controlling its metric, that is, the size and shape of its cells. In this section, we
show that controlling the metric does not control the rank. In the framework of~\cite{HuangRussell2011}, the target metric is
a field $\mathcal M$ called the \emph{monitor}, and a mesh is uniform with respect to it when its cells have equal size in
$\mathcal M$ (\emph{equidistribution}) and are round in $\mathcal M$ (\emph{alignment}). Mesh generation also controls the
condition number $\mathrm{cond}(DX)$ of the chart.

Alignment with a monitor built from $\nabla\phi$ is the natural candidate for a criterion which also controls rank, since it
asks the cells to be shaped by the level sets of $\phi$. However, alignment fixes the shape of the cells but not their
orientation: it determines the Jacobian only up to an orthogonal factor (Proposition~\ref{prop:alignconformal}), and
that factor decides the rank. Theorem~\ref{thm:counterexample} gives a chart which is exactly aligned and equidistributed, and on which a
function of $\phi$ has infinite rank. The rotated Gaussian shows that the condition number does not determine the rank
either. What matters instead is foliation compatibility (Definition~\ref{def:foliationcompat}).

\subsection{The framework}

Following~\cite{HuangRussell2011}, let $x=x(\xi)$ map the computational domain $\Omega_c$ to $\Omega$, with Jacobian
$J=\partial x/\partial\xi$, and let $\mathcal{M}(x)$ be a symmetric positive-definite matrix field, the monitor. The mesh is
uniform with respect to $\mathcal{M}$ if and only if
\begin{align}
\sqrt{\det\mathcal{M}}\,\det J&=\text{const}
&&\text{(equidistribution)},\label{note:eq:equi}\\[1mm]
\tfrac1d\tr\bigl(J^{-1}\mathcal{M}^{-1}J^{-\mathsf T}\bigr)
&=\bigl[\det\bigl(J^{-1}\mathcal{M}^{-1}J^{-\mathsf T}\bigr)\bigr]^{1/d}
&&\text{(alignment)}.\label{note:eq:align}
\end{align}

The monitor which we use measures lengths across the level sets of $\phi$ with weight $\alpha$ and along them with
weight $\beta$.

\begin{definition}[Separating monitor]\label{note:def:monitor}\label{def:sepmonitor}
For $\phi$ with $\nabla\phi\ne0$ and for $\alpha,\beta>0$, we set
\begin{equation}\label{note:eq:monitor}
\mathcal{M}_\phi
=\alpha\,\frac{\nabla\phi\otimes\nabla\phi}{|\nabla\phi|^2}
+\beta\Bigl(I-\frac{\nabla\phi\otimes\nabla\phi}{|\nabla\phi|^2}\Bigr),
\end{equation}
\end{definition}

\begin{proposition}[{\cite[\S2.1]{Huang2006}}]\label{note:prop:conf}\label{prop:alignconformal}
The alignment condition \eqref{note:eq:align} holds if and only if
\begin{equation}\label{note:eq:conformal}
J^{\mathsf T}\mathcal MJ=\chi\,I
\end{equation}
for a positive scalar field $\chi$. That is, unit squares of the computational grid become round cells in $\mathcal M$.
\end{proposition}

In particular, alignment determines $J$ only up to an orthogonal factor: $J=\sqrt\chi\,\mathcal M^{-1/2}Q$ with $Q$
orthogonal. %

\subsection{The counterexample}\label{sec:counterexample}

The point is simple: a rotated grid is a perfect mesh, but a ridge oblique to the grid is not low rank.

Take $\phi(x)=x_2$, whose level sets are horizontal lines, with constant weights $\alpha,\beta$, so that the monitor
$\mathcal M_\phi=\operatorname{diag}(\beta,\alpha)$ is constant. Every linear chart $X(\xi)=J\xi$ whose cells are equal and
round in $\mathcal M_\phi$ is exactly aligned and equidistributed. Such charts are
$J=\sqrt\chi\,\mathcal M_\phi^{-1/2}Q$, with $\chi$ constant and $Q$ orthogonal. In such a chart, $U=G\circ\phi$ becomes
\[
U\circ X=G(a\xi_1+b\xi_2),
\]
a ridge which is oblique to the computational grid unless $Q$ maps each coordinate axis onto a coordinate axis. An oblique ridge is
not a finite sum of products:

\begin{lemma}[Oblique ridges have infinite rank]\label{note:lem:ridge}\label{lem:obliqueridge}
Let $G$ be real-analytic on $\RR$ and not an exponential polynomial, that is, not a finite sum $\sum_ip_i(z)e^{\lambda_iz}$
with polynomials $p_i$ and $\lambda_i\in\mathbb{C}$. Let $a,b\ne0$, and let $I_1,I_2$ be compact nondegenerate intervals. Then the kernel $f(s,t)=G(as+bt)$ on $I_1\times I_2$
has infinite rank, and $r_\varepsilon(f)\to\infty$ as $\varepsilon\to0$.
\end{lemma}

Appendix~\ref{app:proof:lem:obliqueridge} proves Lemma~\ref{lem:obliqueridge}: a kernel of finite rank forces $G$ to
satisfy a linear ordinary differential equation with constant coefficients. The lemma is a local form of Levi-Civita's
theorem on the equation $f(x+y)=\sum_{k=1}^{n}u_k(x)v_k(y)$~\cite{LeviCivita1913}, and of its form for
translation-invariant spaces~\cite{AnseloneKorevaar1964}. For the Gaussian $G(z)=e^{-z^{2}/w^{2}}$ it can be seen directly:
$G(as+bt)=e^{-a^{2}s^{2}/w^{2}}e^{-b^{2}t^{2}/w^{2}}e^{cst}$ with $c=-2ab/w^{2}\ne0$, and $e^{cst}$ has infinite rank,
because exponentials with distinct rates are linearly independent.

\begin{theorem}[Aligned and equidistributed, yet of infinite rank]\label{note:thm:counter}\label{thm:counterexample}
There exist a parallelogram $\Omega\subset\RR^{2}$, a field $\phi$, a monitor $\mathcal{M}_\phi$ of the form
\eqref{note:eq:monitor}, and a linear diffeomorphism $X$ of the computational square $\Omega_c=[0,1]^{2}$ onto $\Omega$
with the following properties:
\begin{itemize}%
\item $X(\partial\Omega_c)=\partial\Omega$;
\item $X$ satisfies equidistribution and alignment \emph{exactly};
\item for $U=G\circ\phi$ with $G$ as in Lemma~\ref{note:lem:ridge}, the function $U\circ X$ has infinite rank:
$r_\varepsilon(U\circ X)\to\infty$ as $\varepsilon\to0$.
\end{itemize}
\end{theorem}

\begin{proof}
Take the construction above, with $Q$ a rotation by an angle $\psi\notin\tfrac\pi2\mathbb Z$, and $\Omega=X(\Omega_c)$. Then
$a,b\ne0$, and Lemma~\ref{lem:obliqueridge} applies.
\end{proof}

Numerically, for every angle $\psi$, both moving-mesh residuals are at rounding level, and the $\varepsilon$-rank is one only
at $\psi=0^\circ$ and $\psi=90^\circ$ (Table~\ref{tab:counter}). Thus, alignment fixes the shape of the cells but not their orientation relative to the
level sets, and the orientation is what decides the rank. Foliation compatibility is a condition on exactly that
orientation.

\subsection{The condition number is the wrong quantity}\label{sec:condnumber}

Mesh generation also controls the condition number $\mathrm{cond}(DX)$ of the chart. The rotated Gaussian, whose singular
values are known exactly, shows that it does not determine the rank. %

\begin{proposition}[Rotated Gaussian; Mehler's formula]\label{note:prop:gauss}\label{prop:rotgauss}
Let $U(x)=\exp(-\tfrac12 x^{\mathsf T}\Sigma^{-1}x)$ on $\RR^2$, and put $\varrho=\Sigma_{12}/\sqrt{\Sigma_{11}\Sigma_{22}}$.
Suppose that $0<|\varrho|<1$. Then the singular values of $U$, as an $L^2$ kernel, form an exact geometric sequence,
\[
\frac{\sigma_{k+1}}{\sigma_k}=z=\frac{|\varrho|}{1+\sqrt{1-\varrho^2}}
\quad\text{for all }k,
\qquad\text{so}\qquad
r_\varepsilon(U)=\Bigl\lceil\frac{\ln(1/\varepsilon)}{\ln(1/z)}\Bigr\rceil\qquad(0<\varepsilon<1).
\]
For $\varrho=0$, the Gaussian is separable, $z=0$, and $r_\varepsilon=1$. Furthermore, for an axis-aligned Gaussian of
variance ratio $\kappa_\Sigma$ rotated by $45^\circ$, this collapses to $z=(\sqrt{\kappa_\Sigma}-1)/(\sqrt{\kappa_\Sigma}+1)$,
with $\sqrt{\kappa_\Sigma}$ the aspect ratio of the level ellipses. This gives
$r_\varepsilon\approx\tfrac12\sqrt{\kappa_\Sigma}\ln(1/\varepsilon)$.
\end{proposition}

The proof is a change of variables in Mehler's formula~\cite{Mehler1866,Watson1933} for the Hermite functions
(Section~\ref{sm:sec:mehler}).

Now take the rotated Gaussian of Figure~\ref{fig:intro}, with $\kappa_\Sigma=25$ at $30^\circ$. The Cartesian chart and the
chart along the principal axes both have $\mathrm{cond}(DX)=1$, yet their $\varepsilon$-ranks at $\varepsilon=10^{-10}$ are
$50$ and $1$. On the other hand, the whitening chart, which maps the level ellipses to circles, has $\mathrm{cond}(DX)=5$ and
still has $\varepsilon$-rank $1$, in every orientation, since $U\circ X=e^{-|\xi|^{2}}$. Thus, the condition number sees the
shape of the cells but not their orientation. %
For a linear chart $X(\xi)=J\xi$, the rank is set by the correlation coefficient $\varrho$ of $U\circ X$, the Gaussian with
covariance $(J^{\mathsf T}\Sigma^{-1}J)^{-1}$ (Proposition~\ref{prop:rotgauss}), and $\varrho$ depends on the shape and the
orientation together. The whitening chart is not foliation-compatible; it has rank one because the profile is an
exponential, as noted after Theorem~\ref{thm:separability}.

\section{The rank cost of misalignment}\label{sec:misalign}\label{note:sec:stability}\label{note:sec:linear}\label{note:sec:dichotomy}

Theorem~\ref{thm:separability} gives rank one when the chart follows the leaves of the function $\phi$ through which
$U=G\circ\phi$ factors. In practice, $\phi$ is not available. The chart is built from an estimate $\hat\phi$, such as a
detected centerline, a fitted level set, or a numerically integrated flow. %
Here we ask what this mismatch costs in $\varepsilon$-rank.

The answer has three parts. First, the amplitude of the mismatch enters only logarithmically
(Theorem~\ref{thm:rankcost}). Second, its rank enters through a binomial count. For polynomial profiles the count is
attained, so no bound linear in that rank holds in general (Proposition~\ref{prop:binomialattained}). Third, what governs
the cost in practice is the structure of the defect in the level variable (Theorem~\ref{thm:productclosed} and
Corollary~\ref{cor:levelanalytic}). Section~\ref{sec:foliation} shows how to read that
structure from derivatives along the flow which builds the chart. %

We first define the leafwise defect, on which the bounds of this section depend. Let $X$ be foliation-compatible for $\hat\phi$, so $\hat\phi\circ X$ depends on $\xi_d$ alone. Let $\phi$ be the function
through which $U=G\circ\phi$ actually factors. We write the true level function in the chart as
\begin{equation}\label{note:eq:defect}
\Psi:=\phi\circ X=\bar\Psi(\xi_d)+\eta(\xi_d,\xi'),
\end{equation}
where $\bar\Psi$ is the leafwise mean of $\Psi$ and $\eta$ is its \emph{leafwise variation}. Then
\[
\widetilde U=U\circ X=G\bigl(\bar\Psi(\xi_d)+\eta(\xi_d,\xi')\bigr),
\]
and $\eta\equiv0$ recovers Theorem~\ref{note:thm:fc}. %
Heuristically, $\eta$ is the normal displacement between the two foliations, measured in $\phi$-units.

The rank of $\eta$ counts the independent ways in which the estimated leaves drift from the true leaves as one moves
across the band. A systematic, level-independent misalignment has $\eta$ of rank one, while a misalignment which
reorganizes itself at every level has high rank. Therefore, two numbers characterize the defect:
\[
\eta_\infty:=\|\eta\|_\infty,\qquad q:=\rank(\eta).
\]

In computations, $\Psi$ is sampled at fiber points $\xi'_i$, $1\le i\le N$, and levels $\xi_{d,j}$, $1\le j\le M$. We use
the same symbols for the samples. That is, $\bar\Psi_j=\bar\Psi(\xi_{d,j})$, $\eta$ is the $N\times M$ matrix
$\eta_{ij}=\eta(\xi_{d,j},\xi'_i)$, $\eta_\infty=\max_{ij}|\eta_{ij}|$, and $q$ is its rank. Furthermore, the sampled field
is $\widetilde U_{ij}=G(\bar\Psi_j+\eta_{ij})$. If $\eta$ is a sum of $q$ separable terms $a_l(\xi')\,b_l(\xi_d)$, then the
matrix has rank at most $q$. We call the vectors $(a_l(\xi'_i))_i$ its fiber-side factors and $(b_l(\xi_{d,j}))_j$ its
level-side factors.

The results below use only that $\bar\Psi$ depends on $\xi_d$ alone, not that it is the leafwise mean. Thus, any split of
$\Psi$ into a function of $\xi_d$ and a remainder can be used, with the remainder in the role of $\eta$. For example,
Proposition~\ref{prop:tiltexact}(i) takes $\bar\Psi(s)=s$. The mean is a convenient default, and
Section~\ref{sec:foliation} takes it with respect to a fixed probability measure on a leaf.

\subsection{The bound}

We bound the $\varepsilon$-rank via the amplitude and the defect rank. The bound has two
steps. First, fix a level $j$ and expand $u\mapsto G(\bar\Psi_j+u)$ in Chebyshev polynomials on
$[-\eta_\infty,\eta_\infty]$. Suppose that $G$ is analytic within a distance $R$ of the real interval on which it is
evaluated. Then the coefficients decay geometrically at a rate set by $R/\eta_\infty$, so a truncation at a degree $m$ of
order $\log(1/\varepsilon)$ meets the tolerance. Second, the truncated expansion is a sum over $p\le m$ of the entrywise powers $\eta^{\circ p}$. Each power appears with
its columns scaled by weights which depend on the level only. Thus, the rank of the truncated expansion is at most the sum
of the ranks of these Hadamard powers, and Lemma~\ref{lem:hadamardpowers} counts them.

\begin{lemma}[Hadamard powers of a rank-$q$ matrix]\label{note:lem:hadamard}\label{lem:hadamardpowers}
Let $\eta=\sum_{l=1}^{q}a_lb_l^{\mathsf T}$ with $a_l\in\RR^N$, $b_l\in\RR^M$. Let $\mathcal F_p$ be the span of all
$p$-fold entrywise products $a_{l_1}\circ\cdots\circ a_{l_p}$, and let $\mathcal L_p$ be the span of all $p$-fold entrywise
products $b_{l_1}\circ\cdots\circ b_{l_p}$, with $\mathcal F_0=\mathcal L_0=\operatorname{span}\{\mathbf 1\}$. Then
\[
\eta^{\circ p}=\sum_{|\alpha|=p}\frac{p!}{\alpha!}\,a^{\alpha}\,(b^{\alpha})^{\mathsf T},
\qquad a^{\alpha}:=\prod_{l}a_l^{\circ\alpha_l},
\]
and, consequently,
\[
\operatorname{rank}\eta^{\circ p}\le\min\{\dim\mathcal F_p,\ \dim\mathcal L_p\}
\le\min\{\binom{q+p-1}{p},N,M\}.
\]
\end{lemma}

\begin{proof}
Entrywise, we have $\eta_{ij}^{p}=\bigl(\sum_la_{il}b_{jl}\bigr)^{p},$ and the multinomial theorem gives the expansion. There are $\binom{q+p-1}{p}$ multi-indices with $|\alpha|=p$.
Furthermore, $a^\alpha\in\mathcal F_p$ and $b^\alpha\in\mathcal L_p$.
\end{proof}

The bound is attained generically: by~\cite[Thm.~4.5]{DammDietrich2023}, there is an open dense subset of the $N\times M$
matrices of rank $q$ on which $\operatorname{rank}\eta^{\circ p}=\min\{\binom{q+p-1}{p},N,M\}$ for every $p$.

\begin{theorem}[Rank cost of an inexact factorization]\label{note:thm:r1prime}\label{thm:rankcost}
Let $\widetilde U_{ij}=G(\bar\Psi_j+\eta_{ij})$, $1\le i\le N$, $1\le j\le M$, where $\bar\Psi_j=\bar\Psi(\xi_{d,j})$,
$\eta\in\RR^{N\times M}$ has rank at most $q$, and $\eta_\infty\ge\max_{ij}|\eta_{ij}|$. Suppose that $G$ is analytic, and bounded
in modulus by $M_G$, on the $R$-neighborhood in $\mathbb C$ of the interval
$[\min_j\bar\Psi_j,\ \max_j\bar\Psi_j]$, and that $\eta_\infty<R$. Put
\begin{equation}\label{note:eq:rho-m}
\begin{gathered}
\rho=\frac{R}{\eta_\infty}+\sqrt{\Bigl(\frac{R}{\eta_\infty}\Bigr)^{2}-1},  %
~~
\varepsilon'=\frac{\varepsilon\,\|\widetilde U\|_F}{\sqrt{NM}},
~~ 
m=\max\Bigl\{0,\Bigl\lceil\frac{\log\bigl(4M_G/((\rho-1)\varepsilon')\bigr)}{\log\rho}\Bigr\rceil\Bigr\}.
\end{gathered}
\end{equation}
For $\eta_\infty=0$, put $m=0$. Then the relative Frobenius $\varepsilon$-rank satisfies
\begin{equation}\label{note:eq:r1prime-crude}
r_\varepsilon(\widetilde U)\ \le\ \binom{q+m}{m},
\end{equation}
and also, for $q_k:=\min\bigl\{\operatorname{rank}B:\ \|T_k(\eta/\eta_\infty)-B\|_F\le\upsilon_k\bigr\},$ $\upsilon_k=\frac{\varepsilon\,\|\widetilde U\|_F}{2(m+1)\,\|c_k\|_\infty},$ 
\begin{equation}\label{note:eq:r1prime-refined}
r_\varepsilon(\widetilde U)\ \le\ \sum_{k=0}^{m}q_k.
\end{equation}
Here $T_k$ is the Chebyshev polynomial, applied entrywise, and $c_k=(c_k(j))_j$ are the Chebyshev coefficients of
$u\mapsto G(\bar\Psi_j+u)$ on $[-\eta_\infty,\eta_\infty]$. Hence,
$r_\varepsilon(\widetilde U)\le\min\bigl\{\binom{q+m}{m},\sum_{k\le m}q_k\bigr\}$.
\end{theorem}

\begin{figure}[tbhp]
\centering
\includegraphics{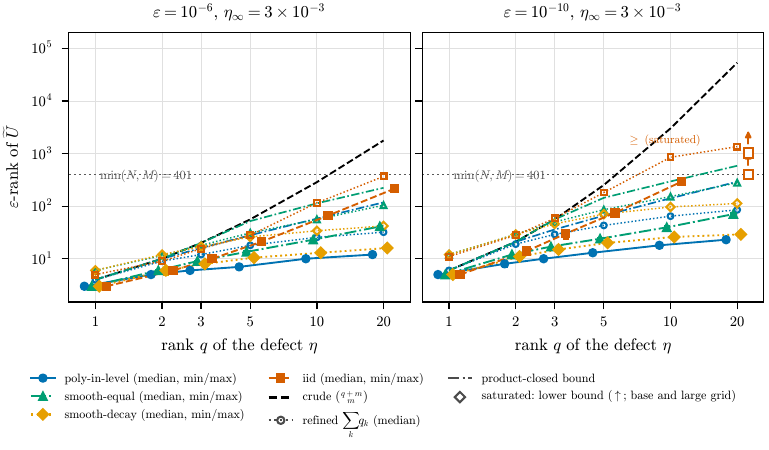}
\caption{$\varepsilon$-rank of $\widetilde U_{ij}=G(\bar\Psi_j+\eta_{ij})$ against the rank $q$ of the defect for four families. Here
$\eta_\infty=3\times10^{-3}$, $G(z)=e^{-z^2/w^2}$, $w=0.12$, with the base grid of $M=401$ levels and $N=1024$ fiber points (where marked,
also the large grid, $M=1024$, $N=2048$), and 10 seeds. Shown are the median and min/max of the observed rank, with the crude bound
$\binom{q+m}{m}$, the refined bound $\sum_kq_k$ (open markers, median), and Theorem~\ref{thm:productclosed} (poly, smooth-equal). Saturated
cells (rank $=\min(N,M)$ on the grid used; iid, $q=20$, $\varepsilon=10^{-10}$) are lower bounds and drawn as open markers with
arrows.}\label{fig:e4}\label{note:tab:r1prime}\label{note:tab:structured}
\end{figure}

The idea is to truncate a Chebyshev expansion in the defect, and then to count the ranks of
the Hadamard powers of $\eta$ which the truncation produces. The full proof is in
Appendix~\ref{app:proof:thm:rankcost}. Neither bound dominates the other; in Figure~\ref{fig:e4}, the refined bound is
the larger at $q=1,2$. The refined bound is computed from the measured defect, so it is an a posteriori bound.

Reading the bound, the amplitude enters only through $\rho$, which is close to $2R/\eta_\infty$ when $\eta_\infty\ll R$. By
\eqref{note:eq:rho-m}, $m$ is $\log(1/\varepsilon')/\log\rho$ up to the additive term $\log\bigl(4M_G/(\rho-1)\bigr)$ in
the numerator. Thus, it grows like $\log(1/\varepsilon')/\log(2R/\eta_\infty)$. Halving $\eta_\infty$ adds about $\log2$
to the denominator. That is, a smaller amplitude lowers the truncation degree, but only slowly. The rank $q$ enters through the count: $\binom{q+m}{m}$ is a polynomial of degree $m$ in $q$. Thus, the
crude bound grows quickly with $q$, even when the amplitude is small. Proposition~\ref{prop:binomialattained} shows that
for general rank-$q$ defects this growth cannot be avoided, and Theorem~\ref{thm:productclosed} shows when it can.

For the Gaussian profile $G(z)=e^{-z^2/w^2}$ of Figure~\ref{fig:e4}, with $w=0.12$, $\eta_\infty=3\times10^{-3}$, $R=w$ and $M_G=e$
(since $|G(x+iy)|=e^{(y^2-x^2)/w^2}\le e$ for $|y|\le w$), and with $401$ equispaced levels $\bar\Psi_j$ on
$[-s_{\max},s_{\max}]$, $s_{\max}=3w/\sqrt2=0.2546$, \eqref{note:eq:rho-m} gives $m=3$ at $\varepsilon=10^{-6}$ and
$m=5$ at $\varepsilon=10^{-10}$. Thus, the crude
bound~\eqref{note:eq:r1prime-crude} is $\binom{q+3}{3}$ and $\binom{q+5}{5}$, respectively.

\subsection{First example: the offset chart}\label{sec:offsetchart}

Here we treat the simplest misaligned chart, an offset chart built on a perturbed centerline. That is, the curve is
displaced along its normal by a small function $\delta$, and the chart follows the normals of the displaced curve. If the
defect were $\delta$ itself (the leading-order model below), it would not depend on the level variable, so $q=1$, and
Theorem~\ref{thm:rankcost} would give $r_\varepsilon\le m+1$. The exact chart of a straight centerline adds one separable
term, so $q=2$ (Proposition~\ref{prop:tiltexact}). That term comes from the tilt of the offset normals.

We work in the planar case for concreteness. %
A centerline is a closed curve in $\RR^2$ or, like the straight and the sinusoidal ones in our examples, an
$L_{\mathrm{per}}$-periodic curve, which we regard as a closed curve on the flat cylinder
$(\RR/L_{\mathrm{per}}\mathbb Z)\times\RR$; reach, offsets, and nearest points are then taken on the cylinder. Let
$\gamma:\TT\to\RR^2$, with $\TT=\RR/L_{\mathrm{per}}\mathbb Z$, be the exact centerline (into the cylinder in the
periodic case), with unit normal $\normal$, let
\begin{equation}\label{note:eq:offset}
X(\theta,s)=\gamma(\theta)+s\,\normal(\theta),\qquad |s|\le s_{\max}<\rch(\gamma),
\end{equation}
and let $U=G\circ d_\gamma$ with $d_\gamma$ the signed distance to $\gamma$, so that $U(X(\theta,s))=G(s)$ exactly. Here
$\rch(\gamma)$ is the reach of $\gamma$, and the hypothesis $s_{\max}<\rch(\gamma)$ makes $\gamma(\theta)$ the nearest
point of $\gamma(\theta)+s\,\normal(\theta)$ on $\gamma$. Let $\hat\gamma=\gamma+\delta\,\normal$ be a perturbed
centerline with normal displacement $\delta:\TT\to\RR$, and let $\widehat X$ be the offset chart built on $\hat\gamma$. We
write $\widehat{\mathbf A}_{ij}=U(\widehat X(\theta_i,s_j))$.

The pulled-back field has the form $U(\widehat X(\theta,s))=G\bigl(s+\eta(\theta,s)\bigr)$, with $\eta$ the
\emph{effective shift}. The \emph{leading-order model} replaces $\eta$ by $\delta(\theta)$. We write
$\widehat{\mathbf A}^{(0)}_{ij}=G\bigl(s_j+\delta(\theta_i)\bigr)$ for its matrix and $A=\|\delta\|_\infty$ for the
amplitude. Proposition~\ref{note:prop:tilt} computes $\eta$ exactly for a straight centerline.

The exact offset chart follows the normals of $\hat\gamma$, which tilt with the slope of $\delta$. This tilt is where the
level variable enters the defect.

\begin{proposition}[Exact chart of a straight centerline]\label{note:prop:tilt}\label{prop:tiltexact}
Let $\gamma$ be the line $\{x_2=0\}$ on the cylinder, parametrized by $\theta\in\TT$ (so that
$d_\gamma=x_2$ and $U=G(x_2)$), let $\delta\in C^1$, and put $\ell=\sqrt{1+\delta'^2}$. The offset chart on the
perturbed centerline $\hat\gamma=(\theta,\delta(\theta))$ is
\[
\widehat X(\theta,s)=\bigl(\theta-s\delta'/\ell,\ \delta+s/\ell\bigr),
\qquad\text{and}
\]
\[
U\bigl(\widehat X(\theta,s)\bigr)=G\bigl(s+\eta(\theta,s)\bigr),\qquad
\eta(\theta,s)=\delta(\theta)-s\,f(\theta),\qquad
f=1-\ell^{-1}\in\bigl[0,\tfrac12\delta'^2\bigr].
\]
Consequently:
\begin{enumerate}[label=(\roman*),leftmargin=*]
\item $\eta$ is a sum of two separable terms, so Theorem~\ref{note:thm:r1prime} applies with $q=2$ and
$\eta_\infty=A+\tau$, where $\tau:=s_{\max}\|f\|_\infty$ (take $\bar\Psi(s)=s$; the proof uses only that $\bar\Psi$
depends on the level variable alone);
\item the singular values of $\widehat{\mathbf A}$ depend on $\delta$ only through the joint distribution of the pairs
$\bigl(\delta(\theta_i),\,\delta'(\theta_i)^2\bigr)$;
\item $\tau\le\tfrac12\,s_{\max}\|\delta'\|_\infty^2$, and, if moreover $\delta''$ is bounded,
$\tau\le A\,s_{\max}\|\delta''\|_\infty$.
\end{enumerate}
\end{proposition}

\begin{proof}
Since $\hat\normal=(-\delta',1)/\ell$, the $x_2$-coordinate of $\widehat X$ is $\delta+s/\ell$, and $U=G(x_2)$. Parts (i)
and (ii) are read off the formula: row $i$ of $\widehat{\mathbf A}$ depends on $(\delta_i,\ell_i)$ alone, and
$\eta(\theta_i,s_j)=\delta_i\cdot1-f_i\,s_j$ has rank two. For (iii), the inequality $1-(1+x)^{-1/2}\le x/2$ for $x\ge0$ gives
the first inequality, and the Landau--Kolmogorov inequality $\|\delta'\|_\infty^2\le2\|\delta\|_\infty\|\delta''\|_\infty$
on $\RR$~\cite{Landau1914,Kolmogorov1939}, hence on $\TT$, gives the second.
\end{proof}

\begin{figure}[tbhp]
\centering
\includegraphics{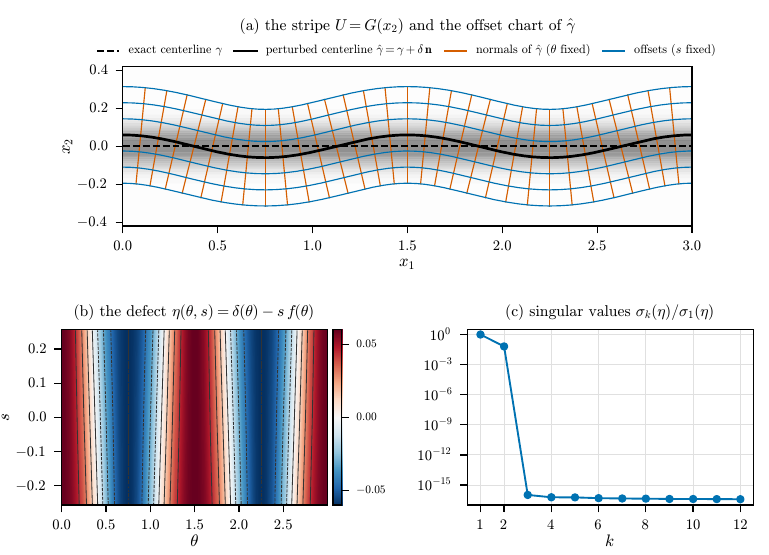}
\caption{The offset chart of a perturbed straight centerline (Proposition~\ref{prop:tiltexact}): the stripe $U=G(x_2)$,
$G(z)=e^{-z^2/w^2}$, $w=0.12$, $s_{\max}=3w/\sqrt2$, period $L_{\mathrm{per}}=3$, perturbation $\delta=A\cos(2\pi K\theta/L_{\mathrm{per}})$ with $A=0.06$ and $K=2$ (the amplitude
is exaggerated for visibility; the experiments use $A\le3\times10^{-3}$). (a) The exact centerline $\gamma$, the perturbed centerline $\hat\gamma$,
and its offset chart; the normals tilt with the slope of $\delta$. (b) The defect $\eta(\theta,s)$. (c) Its singular values: two separable terms,
then rounding error.}\label{fig:misalign}
\end{figure}

Figure~\ref{fig:misalign} shows the
chart, the defect, and its singular values: two separable terms, then rounding error.

\subsection{The binomial count is attained}

For fixed $m$, the crude bound~\eqref{note:eq:r1prime-crude} grows like $q^m/m!$. It is natural to ask whether a bound linear in $q$ holds for every rank-$q$ defect. It does not. %

\begin{proposition}[The binomial bound is attained]\label{note:prop:sharp}\label{prop:binomialattained}
Let $G$ be a polynomial of degree $m$, and let $N,M\ge\binom{q+m}{m}$. Consider generic $a_i,b_j\in\RR^q$ and
$\bar\Psi_j\in\RR$. Then the matrix $\widetilde U_{ij}=G(\bar\Psi_j+a_i\!\cdot\!b_j)$, whose defect $\eta_{ij}=a_i\cdot b_j$
has rank $q$, has rank exactly $\binom{q+m}{m}$. In particular, no bound $r\le c(m)\,q$ holds for all rank-$q$ defects once
$m\ge2$.
\end{proposition}

\begin{proof}
We Taylor-expand $G(\bar\Psi_j+u)$ in $u$ and apply the multinomial theorem to
$u=a_i\cdot b_j$. This gives
\[
\widetilde U_{ij}=\sum_{|\alpha|\le m}a_i^{\alpha}\,c_\alpha(j)
\]
with $c_\alpha(j)=b_j^{\alpha}\,G^{(|\alpha|)}(\bar\Psi_j)/\alpha!$, so $\operatorname{rank}\widetilde U\le\binom{q+m}{m}$.

Monomials are linearly independent functions. Thus, for generic $a_1,\dots,a_N$ with
$N\ge\binom{q+m}{m}$, the vectors $(a_i^\alpha)_i$, $|\alpha|\le m$, are linearly independent. Suppose that $\sum_\alpha\lambda_\alpha c_\alpha(j)=0$ for all $j$. The left side is
the value at $(b_j,\bar\Psi_j)$ of 
\[F=\sum_\alpha\lambda_\alpha f_\alpha, \text{ where } f_\alpha(b,\bar\Psi):=b^\alpha G^{(|\alpha|)}(\bar\Psi)/\alpha!.\] Here $F$ is an element of the fixed space
$W=\operatorname{span}\{f_\alpha:|\alpha|\le m\}$ of dimension $D\le\binom{q+m}{m}$. Consider points $(b_j,\bar\Psi_j)$ in general position, that is, outside a proper analytic subset of the parameters, and let $M\ge D$. Then the evaluation map
$W\to\RR^M$ is injective, because the evaluation functionals at general points span $W^{*}$. Hence, $F\equiv0$.

Next, the coefficient of $b^\alpha$ in $F$ is \[\lambda_\alpha G^{(|\alpha|)}(\bar\Psi)/\alpha!.\] Furthermore,
$G^{(|\alpha|)}\not\equiv0$ because $|\alpha|\le m=\deg G$. Hence, $\lambda=0$, the $c_\alpha$ are independent too, and
$\operatorname{rank}\widetilde U=\binom{q+m}{m}$.
\end{proof}

For $\bar\Psi\equiv0$, Proposition~\ref{note:prop:sharp} is known: the generic rank of $\sum_{k\le m}c_k\eta^{\circ k}$ with
all $c_k\ne0$ is $\min\{\binom{q+m}{m},N,M\}$~\cite[Thm.~3.2]{MaddenThrampoulidis2024}, and that of $\eta^{\circ m}$ is
$\min\{\binom{q+m-1}{m},N,M\}$~\cite[Thm.~4.5]{DammDietrich2023}. Generic shifts $\bar\Psi_j$ make every power up to $m$
contribute, so that only $\deg G=m$ is needed.
Proposition~\ref{note:prop:sharp} is an exact statement about polynomial $G$. For the Gaussian profile, the corresponding
statement about the $\varepsilon$-rank is empirical, and Figure~\ref{fig:e4} reports it.

\subsection{Structure, not rank, decides the cost}

Proposition~\ref{prop:binomialattained} rules out a bound linear in $q$ alone. Lemma~\ref{lem:hadamardpowers} shows what
can replace it: the rank of $\eta^{\circ p}$ is at most the dimension of the span of the $p$-fold products of the level-side
factors, and the binomial count is what remains when nothing is known about them. We therefore place a hypothesis on the
level-side factors $b_l$, which describe how the defect changes from one leaf to the next, and none on the fiber-side
factors $a_l$, which describe how it varies along a leaf and can be as rough as the estimate itself.

\begin{theorem}[Linear rank for product-closed defects]\label{note:thm:structured}\label{thm:productclosed}
Assume the setting of Theorem~\ref{note:thm:r1prime}. Suppose that the level-side factors $b_1,\dots,b_q\in\RR^M$ of
$\eta$ lie in a subspace $\mathcal V\subset\RR^M$. Let $\mathcal V_p$ be the span of the $p$-fold entrywise products of
elements of $\mathcal V$, with $\mathcal V_0=\operatorname{span}\{\mathbf 1\}$, and let $g(p):=\dim\mathcal V_p$. Then
\begin{equation}\label{note:eq:structured}
r_\varepsilon(\widetilde U)\ \le\ \sum_{p=0}^{m}\min\Bigl\{g(p),\ \binom{q+p-1}{p}\Bigr\}.
\end{equation}
The same holds with $\mathcal V\subset\RR^N$ containing the fiber-side factors. In particular:
\begin{enumerate}[label=(\alph*),leftmargin=*]
\item if every $b_l$ samples a polynomial of degree $<n$ in the level variable, then $g(p)\le p(n-1)+1$ and
$r_\varepsilon(\widetilde U)\le(m+1)\bigl(1+\tfrac{(n-1)m}{2}\bigr)$;
\item if every $b_l$ samples a trigonometric polynomial of degree $\le n$, then $g(p)\le2pn+1$ and
$r_\varepsilon(\widetilde U)\le(m+1)(1+nm)$.
\end{enumerate}
In both cases, the bound is linear in the dimension of the factor space. Furthermore, no hypothesis is made on the
fiber-side factors $a_l$, which may be as rough as one likes.
\end{theorem}

\begin{proof}
By
Lemma~\ref{note:lem:hadamard}, $\operatorname{rank}\eta^{\circ p}\le\min\{\dim\mathcal L_p,\binom{q+p-1}{p}\}$.
Furthermore, $\mathcal L_p\subset\mathcal V_p$. We insert this in the argument in the proof of Theorem~\ref{note:thm:r1prime}. For the dimensions $g(p)$, a product of $p$ polynomials of degree $<n$ has degree $\le p(n-1)$. Similarly, a product of
$p$ trigonometric polynomials of degree $\le n$ has degree $\le pn$. These give the bounds on $g(p)$. Finally,
\[\sum_{p=0}^{m}\bigl(p(n-1)+1\bigr)=(m+1)\bigl(1+(n-1)m/2\bigr) \text{ and } \sum_{p=0}^{m}(2pn+1)=(m+1)(1+nm).\]
\end{proof}

Exactly polynomial dependence on the level variable is special. Suppose that the level-side dependence is only analytic.
Then it can be truncated to a polynomial whose degree grows like $\log(1/\varepsilon)$, and
Theorem~\ref{thm:productclosed}(a) applies to the truncation. Since $m$ also grows like $\log(1/\varepsilon)$, the rank
stays polylogarithmic in $1/\varepsilon$, still with no assumption on the fiber side.

\begin{corollary}[Level-analytic defects: polylogarithmic rank]\label{note:cor:analytic}\label{cor:levelanalytic}
Let $\eta_{ij}=\eta_i(x_j)$ with level variable $x_j\in[-1,1]$. Suppose that every $\eta_i$ extends analytically to one
and the same Bernstein ellipse $E_{\rho_s}$, that is, the open elliptical region with foci $\pm1$ and semi-axis sum
$\rho_s$. Suppose also that $|\eta_i|\le M_\eta$ there, uniformly in $i$, with no assumption on the dependence on $i$.

Let $\upsilon>0$, and let $M_1$ bound $|G'|$ on the real interval
$[\min_j\bar\Psi_j-\eta_\infty-\upsilon,\ \max_j\bar\Psi_j+\eta_\infty+\upsilon]$, with
$4M_1\upsilon\sqrt{NM}\le\varepsilon\|\widetilde U\|_F$. Thus, $\upsilon\approx\varepsilon\bar U/(4M_1)$, where
$\bar U=\|\widetilde U\|_F/\sqrt{NM}$. Assume the setting of Theorem~\ref{thm:rankcost}, with $\eta_\infty+\upsilon<R$, and put
\[
n_\upsilon=\max\Bigl\{1,\ 1+\Bigl\lceil\frac{\log\bigl(2M_\eta/((\rho_s-1)\upsilon)\bigr)}{\log\rho_s}\Bigr\rceil\Bigr\}.
\]
Let $\eta_\upsilon$ be the $N\times M$ matrix whose entries $\eta_{\upsilon,ij}$ are the values at $x_j$ of the Chebyshev
series of $\eta_i$ truncated at degree $n_\upsilon-1$. Then, for $\varepsilon\le1$,
\[
r_\varepsilon(\widetilde U)\le(m+1)\bigl(1+\tfrac{(n_\upsilon-1)m}{2}\bigr)=O\bigl(\log^{3}(1/\varepsilon)\bigr),
\]
where $m$ is computed from \eqref{note:eq:rho-m} with $(\varepsilon/4,\ \eta_\infty+\upsilon)$ in place of
$(\varepsilon,\ \eta_\infty)$ and with $\varepsilon'$ formed from $\|G(\bar\Psi+\eta_\upsilon)\|_F$.
\end{corollary}

The proof is in Appendix~\ref{app:proof:cor:levelanalytic}. It truncates the Chebyshev series of each $\eta_i$ at degree
$n_\upsilon-1$, with error at most $\upsilon$ by~\cite[Thm.~8.1]{Trefethen2019}. Then it bounds the resulting change of
$\widetilde U$ by the mean value theorem and applies Theorem~\ref{thm:productclosed}(a) to the truncated defect
at tolerance $\varepsilon/4$.

The question left open is when the defect of a chart has simple level-side structure, without computing the defect.
Section~\ref{sec:foliation} answers it through derivatives of $\phi$ along the flow which builds the chart.

\section{Where the cost comes from: the foliation order}\label{sec:foliation}\label{note:sec:foliation}

By Section~\ref{sec:misalign}, the cost of a defect is governed by its structure in the level variable
(Theorem~\ref{thm:productclosed} and Corollary~\ref{cor:levelanalytic}). Here, we show that this structure is
controlled by the derivatives of $\phi$ along the flow which builds the chart, without computing $\eta$.

Let $\hat\phi$ and $\phi$ be smooth with nonvanishing gradients on the band, and let
$\hat V:=\nabla\hat\phi/|\nabla\hat\phi|^{2}$, so that $\hat V\hat\phi=1$. Let
\[X(\xi_d,y)=\Phi^{\hat V}_{\xi_d-c_\ast}(y),~~ y\in\hat F=\hat\phi^{-1}(c_\ast),\] be the chart built by the flow of
$\hat V$; see Theorem~\ref{thm:bandproduct}. Thus, $\hat\phi\circ X=\xi_d$. We write
$\phi\circ X=\bar\Psi(\xi_d)+\eta(\xi_d,y)$ as in Section~\ref{note:sec:stability}, where $\bar\Psi$ is the leafwise mean
with respect to a fixed probability measure $\mathbb P$ on $\hat F$. We write $\hat V^{n}\phi$ for the
$n$-th derivative of $\phi$ along the flow lines, that is, $\hat V$ applied $n$ times.

Since the chart is the flow of $\hat V$, we have $\partial_{\xi_d}(f\circ X)=(\hat Vf)\circ X$ for every smooth $f$: differentiating
in the level coordinate is differentiating along the flow lines. Therefore, the dependence of $\eta$ on the level coordinate
is carried by the derivatives $\hat V^{n}\phi$, which can be computed without forming $\eta$.

\begin{proposition}[Derivatives of the defect along the flow]\label{note:prop:lie}\label{prop:liederivs}
The following hold:
\begin{enumerate}[label=(\roman*),leftmargin=*]
\item for $n\ge1$, $\partial_{\xi_d}^{n}(\phi\circ X)=(\hat V^{n}\phi)\circ X$; in particular,
$\partial_{\xi_d}(\phi\circ X)=\lambda\circ X$ with
\[
\lambda:=\hat V\phi=\frac{\langle\nabla\phi,\nabla\hat\phi\rangle}{|\nabla\hat\phi|^{2}}
=\frac{|\nabla\phi|}{|\nabla\hat\phi|}\cos\zeta ,
\]
where $\zeta$ is the angle between the true and the estimated gradient;
\item $\partial_{\xi_d}^{n}\eta=(\hat V^{n}\phi)\circ X-\bar\Psi^{(n)}(\xi_d)$, where $\bar\Psi^{(n)}(\xi_d)$ is the
$\mathbb P$-mean over $\hat F$ of the first term;
\item for $k\ge1$, the following are equivalent: (a) for every $y$, $\xi_d\mapsto\eta(\xi_d,y)$ is a polynomial of degree
$<k$; (b) $\hat V^{k}\phi$ is constant on every leaf of $\hat\phi$, that is, it is a function of $\hat\phi$ alone.
\end{enumerate}
\end{proposition}

\begin{proof}
The idea is that the level coordinate of the chart is the time of the flow of $\hat V$.

(i) Since $X$ is the flow of $\hat V$ in time $\xi_d-c_\ast$, we have $\partial_{\xi_d}(f\circ X)=(\hat Vf)\circ X$ for
every smooth $f$. Iterating gives the first formula. The formulas for $\lambda$ are the definition of $\hat V$.

(ii) We subtract $\bar\Psi^{(n)}$, differentiating $\bar\Psi(\xi_d)=\int\phi\circ X(\xi_d,y)\,d\mathbb P(y)$ under the
integral.

(iii) On an interval, (a) holds if and only if $\partial_{\xi_d}^{k}\eta(\cdot,y)\equiv0$ for each $y$. By (ii), this holds
if and only if $(\hat V^{k}\phi)\circ X(\xi_d,y)=\bar\Psi^{(k)}(\xi_d)$ does not depend on $y$. Since $X(\xi_d,\cdot)$
maps $\hat F$ onto the leaf $\hat\phi^{-1}(\xi_d)$, this says that $\hat V^k\phi$ is constant on each leaf.
\end{proof}

\begin{definition}[Foliation order]\label{def:folorder}
The \emph{foliation order} $\folord(\phi\,|\,\hat\phi)$ of $\phi$ relative to $\hat\phi$ is the least $k\ge1$ for which
the equivalent conditions of Proposition~\ref{prop:liederivs}(iii) hold, and $\infty$ if there is none.
\end{definition}

When the foliation order is finite, Theorem~\ref{note:thm:structured}(a) with $n=\folord(\phi\,|\,\hat\phi)$ gives
\[
r_\varepsilon(\widetilde U)\le(m+1)\Bigl(1+\frac{(\folord(\phi\,|\,\hat\phi)-1)m}{2}\Bigr),
\]
with no assumption on how rough the foliations are in the fiber direction. The first two orders have a direct meaning.
Order $1$ means that the gradient ratio $\lambda$ is a function of $\hat\phi$ alone. That is, the misalignment is a pure
displacement of the leaves, $\eta=\eta(y)$, which is the leading-order model of Section~\ref{note:sec:dichotomy}. Order $2$
means that $\hat V\lambda$ is a function of $\hat\phi$ alone.

\begin{corollary}[The offset chart of a straight centerline has order $2$]\label{cor:orderstraight}
Assume the setting of Proposition~\ref{prop:tiltexact}, with $\delta\in C^{2}$ and $s_{\max}<\rch(\hat\gamma)$. Let $\phi=x_2$, and let
$\hat\phi$ be the signed distance to $\hat\gamma$. Then $\folord(\phi\,|\,\hat\phi)=2$, unless $\delta$ is constant on $\TT$,
in which case it is $1$.
\end{corollary}

\begin{proof}
Inside the reach, the nearest point of $\widehat X(\theta,s)$ on $\hat\gamma$ is
$\hat\gamma(\theta)$~\cite[Thm.~4.8]{Federer1959} (on the cylinder, applied to the periodic lift of $\hat\gamma$ in
$\RR^2$, as after Proposition~\ref{prop:injectivity}). Thus, $\hat\phi\circ\widehat X=s$ and $|\nabla\hat\phi|=1$.
Moreover, $\widehat X$ is a $C^1$ diffeomorphism onto its image (Proposition~\ref{prop:injectivity}(iii) below) and
$\nabla\hat\phi=\hat\normal\circ\theta\circ\widehat X^{-1}$, so $\hat\phi\in C^2$ there, which is all that
Proposition~\ref{prop:liederivs} uses for $n\le2$.
Furthermore, $\hat V=\nabla\hat\phi=\hat\normal(\text{foot point})$ is constant along each normal line. Hence,
$\lambda=\hat V\phi=\langle e_2,\hat\normal\rangle=1/\ell(\theta)$ is constant along $\hat V$. Therefore,
$\hat V^{2}\phi=\hat V\lambda=0$ is constant on every leaf, and the order is at most $2$ by
Proposition~\ref{prop:liederivs}(iii).

It remains to decide when the order is $1$. This is the case if and only if $\lambda=1/\ell$ is constant on every leaf of
$\hat\phi$, that is, if and only if $\ell=\sqrt{1+\delta'^2}$ is constant. Equivalently, $|\delta'|$ is constant. Since
$\delta'$ is continuous and periodic with mean zero, this holds if and only if $\delta'=0$. %
\end{proof}

\begin{remark}
Corollary~\ref{cor:orderstraight} is the geometric reason behind Proposition~\ref{prop:tiltexact}. Indeed, by
Proposition~\ref{prop:liederivs}(iii), order two means that $\eta(\cdot,y)$ is affine in the level coordinate for every
$y$. Thus, whatever $\delta$ is, $\eta$ is a sum of two separable terms, one constant and one linear in the level
coordinate. Up to the leafwise mean, these are the displacement $\delta$ and the tilt $-s\,f$ of
Proposition~\ref{prop:tiltexact}. The frequency content of $\delta$ enters only through the fiber-side factors, on which
Theorem~\ref{thm:productclosed} makes no assumption.
\end{remark}

\begin{remark}[Reparametrizing the level function]\label{rem:folorderinvariance}
The foliation order depends on $\hat\phi$ as a function, not only on its leaves. If $\beta$ is smooth and strictly
monotone, then $\beta\circ\hat\phi$ has the same leaves, and its flow follows the same curves at a different speed. On the
same set of leaves, the chart therefore samples the same points, and the sampled matrix is unchanged. Its order, however,
can change: for the offset chart of Corollary~\ref{cor:orderstraight}, replacing the signed distance $s$ by $\beta(s)=s+s^3$
makes the order infinite, unless $\delta$ is constant: in the coordinate $t=\beta(s)$, the defect is affine in
$\beta^{-1}(t)$, which is not a polynomial. Since the rank bound depends only on the sampled matrix, it holds with the least order over all
such $\beta$, which depends on the two foliations alone. Similarly, $U=G\circ\phi=(G\circ\alpha^{-1})\circ(\alpha\circ\phi)$
for every strictly monotone $\alpha$, and one may use the representative of $\phi$ with the smallest order, at the price of
changing $R$ and $M_G$ in Theorem~\ref{thm:rankcost}. For a monotone profile, $\alpha=G$ gives $U$ itself, so the order can
then be computed from the field. %
\end{remark}
\begin{remark}[The gauge is free only when the factorization is exact]\label{note:rem:gauge-defect}\label{rem:gaugedefect}
Proposition~\ref{prop:gaugefree} is a statement about $\eta=0$. Suppose that $U=G\circ\phi$, but the chart follows the
level sets of an estimate $\hat\phi$. Then a change of fiber coordinates $y\mapsto\pi_{\xi_d}(y)$ replaces $\eta(\xi_d,y)$
by $\eta(\xi_d,\pi_{\xi_d}(y))$. Thus, it changes the dependence of the defect on the level variable, which is what the
rank sees (Corollary~\ref{note:cor:analytic}).

For example, a defect $\eta=A\sin Ky$ is of rank one and order one in the flow-transported gauge. In the gauge
$y=\xi'+b\,\xi_d\sin\xi'$, with $\xi_d$ normalized to $[-1,1]$, its level dependence has up to
$O\bigl(Kb+\log(1/\varepsilon)\bigr)$ significant modes, by the Jacobi--Anger expansion. Numerically, for $K=64$, the $\varepsilon$-rank of $\widetilde U$
rises from $4$ at $b=0$ to $120$ at $b=0.9$. For $\eta=0$, the
$\varepsilon$-rank is $1$ in every gauge. Thus, a gauge may be chosen freely only when $\eta=0$. Otherwise, it should be
flow-transported, for which the level dependence is that of Proposition~\ref{note:prop:lie}, or it should at least keep
$\eta(\xi_d,\pi_{\xi_d}(y))$ simple in $\xi_d$.
\end{remark}

For a curved centerline, $\lambda$ in general varies along the normals, and the order can be infinite, for instance for a
circle. The dependence on the level coordinate is still analytic, and the relevant
quantity is then its radius of analyticity. In words, the next result says: if the derivatives of $\phi$ along the flow
vary little across each leaf, and grow no faster than those of an analytic function with radius $R_\phi$, then the defect
is analytic in the level coordinate, and Corollary~\ref{cor:levelanalytic} gives a rank polylogarithmic in $1/\varepsilon$.

\begin{corollary}[A geometric bound on the rank]\label{note:cor:geometric}\label{cor:georank}
Let the level interval of the band be $[c-L,c+L]$. Suppose that there are $C>0$ and $R_\phi>L$ such that the derivatives
along the flow satisfy $\operatorname{osc}_{\text{leaf}}\hat V^{n}\phi\le C\,n!\,R_\phi^{-n}$ for all $n\ge1$ \emph{on every
leaf of the band}, where $\operatorname{osc}_{\text{leaf}}f$ is the difference between the largest and the smallest value of
$f$ on a leaf. Suppose also that $\operatorname{osc}\,\phi$ on the leaf at level $c$ is
at most $\eta_c$. Then the analyticity hypothesis of Corollary~\ref{note:cor:analytic} holds with
\[
\rho_s=\frac{R'_\phi}{L}+\sqrt{\Bigl(\frac{R'_\phi}{L}\Bigr)^{2}-1},\qquad
M_\eta=\eta_c+\frac{C\,R'_\phi}{R_\phi-R'_\phi},\qquad\text{for any }L<R'_\phi<R_\phi,
\]
so, under the other hypotheses of that corollary (on $G$, and $\eta_\infty+\upsilon<R$),
$r_\varepsilon(\widetilde U)=O\bigl(\log^{3}(1/\varepsilon)\bigr)$, and the dimension of the defect is
$n_\upsilon\simeq\log(M_\eta/\upsilon)/\log\rho_s$. Furthermore, suppose that $\eta(\cdot,y)$ is known to be
real-analytic on $[c-L,c+L]$ for every $y$, for instance because $\phi$ and $\hat\phi$ are real-analytic and $\hat F$ is
compact. Then the oscillation bound is needed on the central leaf only.
\end{corollary}

\begin{proof}[Proof of Corollary~\ref{cor:georank}]
The idea is to bound the Taylor coefficients of $\eta(\cdot,y)$ at $c$ by the oscillations of $\hat V^{n}\phi$, and then
to fit a Bernstein ellipse inside the disk of convergence.

\emph{Step 1: the derivatives.} By Proposition~\ref{note:prop:lie}(ii),
$\partial^{n}\eta(\xi,y)=(\hat V^{n}\phi)\circ X(\xi,y)-\bar\Psi^{(n)}(\xi)$ is a value minus a $\mathbb P$-mean of values on
the same leaf. Thus, $|\partial^{n}\eta(\xi,y)|\le\operatorname{osc}_{\text{leaf at }\xi}\hat V^{n}\phi$.

\emph{Step 2: the Taylor series.} There are two cases. Suppose first that the oscillation bound holds on every leaf. Then,
for $\xi$ in the interval (where $|\xi-c|\le L<R_\phi$), Taylor's formula with remainder gives
\[
\bigl|\eta(\xi,y)-\sum_{n<N}\partial^{n}\eta(c,y)(\xi-c)^{n}/n!\bigr|\le C(|\xi-c|/R_\phi)^{N}\to0.
\]
Hence, $\eta(\cdot,y)$ is the sum of its Taylor series at $c$.

Suppose instead that $\eta(\cdot,y)$ is real-analytic on the interval. Then it agrees near $c$ with the sum of that series, hence on the whole interval by the identity theorem. That
series converges in $|\xi-c|<R_\phi$ because $|\partial^{n}\eta(c,y)|\le Cn!\,R_\phi^{-n}$.

\emph{Step 3: the extension.} In both cases,
\[
\eta(\xi,y)=\eta(c,y)+\sum_{n\ge1}\partial^{n}\eta(c,y)\,\frac{(\xi-c)^{n}}{n!},
\qquad
|\partial^{n}\eta(c,y)|\le C\,n!\,R_\phi^{-n},\quad|\eta(c,y)|\le\eta_c,
\]
and the series extends $\eta(\cdot,y)$ analytically to $|\xi-c|<R_\phi$. There,
$|\eta|\le\eta_c+C|\xi-c|/(R_\phi-|\xi-c|)$, uniformly in $y$. Finally, the open Bernstein ellipse of $[c-L,c+L]$ with
parameter $\rho_s$ has semi-major axis $R'_\phi$. Thus, it lies in $|\xi-c|<R'_\phi$, where $|\eta|\le M_\eta$ uniformly
in $y$.
\end{proof}

\begin{remark}[Consequence for the design of a chart]
What the rank sees is how the estimation error varies across the leaves, not how rough it is along them. For a chart built
from a field estimate (Section~\ref{sec:invert}), this means that $\lambda=\hat V\phi$ should be nearly constant along each
flow line, with derivatives as in Corollary~\ref{cor:georank} for a radius large compared with the band half-width, and
that the fiber coordinates should be carried along the flow (Remark~\ref{rem:gaugedefect}).
\end{remark}
\section{Invertibility and the design principle}\label{sec:invert}\label{note:sec:inject}

Section~\ref{sec:misalign} measured what a perturbed chart costs in rank. That cost exists whether or not the chart folds,
since the bounds there are statements about the sampled matrix, which is defined for any chart. However, the
reconstruction~\eqref{note:eq:pullback} needs the chart to be invertible. In this section, we ask what the perturbation
does to invertibility, and how a chart can be made invertible by construction.

We first characterize when an offset chart is injective.

\begin{proposition}[Chart injectivity]\label{note:prop:inject}\label{prop:injectivity}
Let $\hat\gamma$ be a closed embedded $C^{2}$ curve in $\RR^{2}$, or on the cylinder of
Section~\ref{sec:offsetchart}, with curvature $\hat\kappa$ with respect to
$\hat\normal$. Let $\widehat X(\theta,s)=\hat\gamma(\theta)+s\,\hat\normal(\theta)$ on $\TT\times[-s_{\max},s_{\max}]$, and
let $\mu:=s_{\max}\|\hat\kappa\|_\infty$. Then the following hold:
\begin{enumerate}[label=(\roman*),leftmargin=*]
\item $\det D\widehat X(\theta,s)=|\hat\gamma'(\theta)|\bigl(1-s\,\hat\kappa(\theta)\bigr)$;
\item $\det D\widehat X\ne0$ on the whole strip if and only if $\mu<1$, and if $\mu>1$, then $\widehat X$ is not
injective;
\item $\widehat X$ is a diffeomorphism onto its image if and only if $s_{\max}<\rch(\hat\gamma)$, where
$\rch(\hat\gamma)=\min\{\omega/2,\ 1/\|\hat\kappa\|_\infty\}$ and $\omega$ is the length of the shortest double normal of
$\hat\gamma$;
\item if $\gamma$ is a straight line and $\hat\gamma=\gamma+\delta\,\normal$, then
$\hat\kappa=\delta''/(1+\delta'^{2})^{3/2}$. Hence, for $\delta(u)=A\cos(2\pi Ku/L_{\mathrm{per}})$ on a period
$L_{\mathrm{per}}$, one has $\|\hat\kappa\|_\infty\ge A(2\pi K/L_{\mathrm{per}})^{2}$, and $\mu\ge1$ as soon as
$s_{\max}A(2\pi K/L_{\mathrm{per}})^{2}\ge1$.
\end{enumerate}
\end{proposition}

Proposition~\ref{prop:injectivity} is classical. Part~(i) is Frenet's formula, and (ii) follows from (i) by invariance of
domain. Part~(iii) combines the formula for the reach~\cite[Thm.~3.4 and Prop.~A.1]{Aamari2019} with the continuity of the
nearest-point map inside the reach~\cite[Thm.~4.8]{Federer1959}. On the cylinder, these results apply to the
$L_{\mathrm{per}}$-periodic lift of $\hat\gamma$ in $\RR^2$: it is a closed set with the same nearest points, curvature,
and double normals, and its periodicity provides the compactness which the formula for the reach uses. Part~(iv) is the curvature $\delta''/(1+\delta'^2)^{3/2}$
of the graph of $\delta$, evaluated where $\delta'=0$. We call $\mu$ the \emph{injectivity margin}. By (iii), the
sharp condition is the reach of $\hat\gamma$~\cite{Federer1959}, which is stronger than a curvature bound: two parts of the
curve can come close to each other without either being curved, and no curvature test detects
this~\cite[Ch.~7]{Breiding2024}.

\subsection{The dichotomy}

In this subsection, we set the rank of an offset chart against its injectivity. The rank half of the following theorem
restates Proposition~\ref{prop:tiltexact} and Theorem~\ref{thm:rankcost}.

\begin{theorem}[Dichotomy]\label{note:thm:dichotomy}\label{thm:rankinjdichotomy}
Let $\delta$ be a normal perturbation of the centerline, let $A=\|\delta\|_\infty$, and let
$\mu=s_{\max}\|\hat\kappa\|_\infty$. Then the perturbation reaches the two halves of the representation through different
quantities and with different failure modes.
\begin{enumerate}[label=(\roman*),leftmargin=*]
\item \emph{The rank degrades softly.} In the exact chart of a straight centerline, the rank bound depends on $\delta$
only through $A$, the tilt amplitude $\tau\le\tfrac12\,s_{\max}\|\delta'\|_\infty^2$, and, through $\varepsilon'$, the
norm $\|\widehat{\mathbf A}\|_F$, with one extra separable term (Proposition~\ref{note:prop:tilt} and
Theorem~\ref{note:thm:r1prime} with $q=2$ and $\eta_\infty=A+\tau$). The amplitudes enter through
$\log\bigl(2R/(A+\tau)\bigr)$, provided $A+\tau<R$. %
\item \emph{The injectivity fails hard.} The Jacobian of the chart vanishes somewhere in the strip exactly when $\mu\ge1$,
and for $\mu>1$ the chart is not injective. In general, the chart is a diffeomorphism if and only if
$s_{\max}<\rch(\hat\gamma)$, and $\rch(\hat\gamma)\le1/\|\hat\kappa\|_\infty$. Thus, a single-valued
pullback~\eqref{note:eq:pullback} is well posed only below this threshold. For a flat centerline with small slope, the
threshold is controlled by $\|\delta''\|_\infty$ (Proposition~\ref{note:prop:inject}).
\end{enumerate}
\end{theorem}

\begin{proof}
(i) By Proposition~\ref{note:prop:tilt}, $\eta=\delta-sf$ has rank at most $2$ and $|\eta|\le A+\tau$, so the crude
bound~\eqref{note:eq:r1prime-crude} with $q=2$ applies with $\rho$ computed from $A+\tau$; Proposition~\ref{note:prop:tilt}(iii)
estimates $\tau$. (ii) This is Proposition~\ref{note:prop:inject}(ii)--(iii).
\end{proof}

\begin{figure}[tbhp]
\centering
\includegraphics{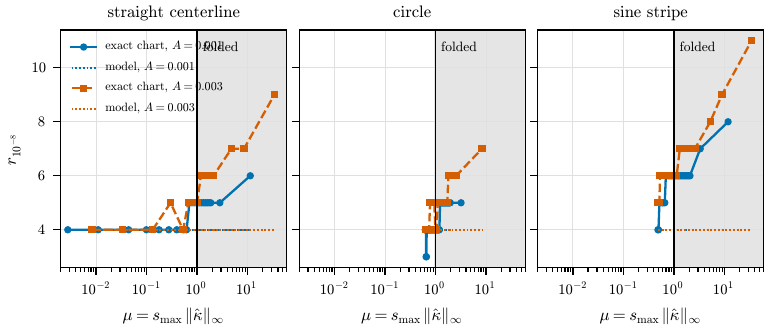}
\caption{Rank against the injectivity margin $\mu=s_{\max}\|\hat\kappa\|_\infty$ for the exact offset chart and the leading-order model
($G(z)=e^{-z^2/w^2}$, $w=0.12$, $s_{\max}=0.6308$, $8192\times401$ nodes, $\delta=A\cos(2\pi K\theta/L_{\mathrm{per}}+0.37)$, $\varepsilon=10^{-8}$). Gray:
the chart folds ($\mu\ge1$).}\label{fig:e3rank}
\end{figure}

Table~\ref{tab:inject} shows three things. First, the model's rank does not move. Second, the exact chart's rank does
move, but slowly: from $4$ to $13$ while $\mu$ rises four-thousandfold on the flat strip. While $\mu<1$, it exceeds the
model's by at most two. Third, the chart ceases to be injective no later than $\mu=1$, whatever the rank does. This is the
dichotomy in quantitative form. Figure~\ref{fig:e3rank} shows the same behavior for the pre-registered sweep of
Section~\ref{sec:numerics}, at a larger offset.

\begin{table}[tbhp]
\centering\scriptsize
\caption{Exact offset chart on $\hat\gamma=\gamma+\delta\,\normal$, $\delta=3\times10^{-3}\cos(2\pi K\theta/L_{\mathrm{per}})$,
$s_{\max}=3w/\sqrt2$, $w=0.12$, $\varepsilon=10^{-8}$, $N_s=401$; flat strip $L_{\mathrm{per}}=1$, $N_\theta=4096$; circle
$R_0=1$, $L_{\mathrm{per}}=2\pi$, $N_\theta=8192$. Here
$\mu=s_{\max}\max|\hat\kappa|$ on the chart's nodes.}\label{tab:inject}\label{note:tab:inject}
\begin{adjustbox}{max width=\textwidth}\begin{tabular}{lrrrlrrr}
\toprule
geometry & $K$ & $\max|\hat\kappa|$ & $\mu$ & label & $r_{\rm model}$ & $r_{\rm exact}$ & $\min\det$ \\
\midrule
flat & 1 & 0.118 & 0.0301 & regular & 4 & 4 & 0.97 \\
flat & 4 & 1.89 & 0.482 & regular & 4 & 5 & 0.518 \\
flat & 5 & 2.96 & 0.754 & regular & 4 & 6 & 0.246 \\
flat & 8 & 7.58 & 1.93 & singular & 4 & 7 & -0.93 \\
flat & 16 & 30.3 & 7.72 & singular & 4 & 8 & -6.72 \\
flat & 64 & 485 & 123 & singular & 4 & 13 & -122 \\
circle & 1 & 1 & 0.255 & regular & 4 & 4 & 0.743 \\
circle & 16 & 1.76 & 0.448 & regular & 4 & 5 & 0.554 \\
circle & 24 & 2.71 & 0.691 & regular & 4 & 5 & 0.31 \\
circle & 32 & 4.05 & 1.03 & singular & 4 & 6 & -0.0312 \\
circle & 64 & 13.2 & 3.36 & singular & 4 & 7 & -2.37 \\
circle & 256 & 197 & 50.1 & singular & 4 & 11 & -49.2 \\
\bottomrule
\end{tabular}
\end{adjustbox}
\end{table}

\subsection{Injectivity reduces to the leaves}

In this subsection, we show that foliation compatibility, which Section~\ref{sec:setting} identified as the condition for
rank one, also settles invertibility. %
Precisely, injectivity of the chart reduces to
injectivity of its leaf charts.

\begin{theorem}[Injectivity reduces to the leaves]\label{note:thm:inject}\label{thm:leafinjectivity}
Let $X$ be foliation-compatible with $\bar\Psi=\phi\circ X$ strictly monotone in $\xi_d$. Suppose that $X$ restricted to
each slice $\{\xi_d=t\}$ is injective into the leaf $\phi^{-1}(\bar\Psi(t))$. Then $X$ is injective on $\Omega_c$. Thus,
injectivity of the $d$-dimensional chart reduces to injectivity of the $(d-1)$-dimensional leaf charts, and no curvature,
reach, or offset-distance hypothesis enters. Furthermore, if $\nabla\phi\ne0$, $\bar\Psi'\ne0$, and every restriction
$X(\cdot,t)$ is an immersion, then $X$ is a diffeomorphism onto its image.
\end{theorem}

\begin{proof}

Suppose that $X(\xi)=X(\tilde\xi)$. Applying $\phi$ gives $\bar\Psi(\xi_d)=\bar\Psi(\tilde\xi_d)$, so $\xi_d=\tilde\xi_d$
by monotonicity. Then injectivity on that slice gives $\xi'=\tilde\xi'$. For the last assertion, the columns $\partial_{\xi_i}X$, $i<d$, span the tangent space of the leaf. On the other hand,
\[
\nabla\phi\cdot\partial_{\xi_d}X=\partial_{\xi_d}(\phi\circ X)=\bar\Psi'(\xi_d)\ne0,
\]
so $\partial_{\xi_d}X$ is not tangent to it. Hence, $DX$ is invertible. Finally, an injective local diffeomorphism is a
diffeomorphism onto its (open) image.
\end{proof}

The hypothesis on the slices is not vacuous. For $d=2$, it asks for a monotone parametrization of each level curve. However, for $d\ge3$, it is a
surface-parametrization problem on each leaf, subject to the topological obstruction discussed in Section~\ref{sec:discussion}.

For offset charts, this gives a new reading of the reach condition. Inside the reach, the offsets of $\hat\gamma$ are
exactly the level sets of the signed distance to $\hat\gamma$, so the
offset chart is foliation-compatible, but for that distance function, not for the field. Thus, the reach condition is not a
limitation of charts adapted to level sets; it is the price of building one by offsetting an estimated curve, rather than
by following the level sets of an estimated field.

\subsection{The construction}

Theorem~\ref{thm:leafinjectivity} and Proposition~\ref{prop:gaugefree} suggest building the chart from an estimate of the
field rather than from a curve. The construction has four steps.
\begin{enumerate}[label=\textbf{C\arabic*.},leftmargin=*]
\item \emph{Estimate the field.} Form a smooth interpolant $\hat\phi$ of the sampled data, such as a tensor-product spline,
or a Chebyshev or Fourier interpolant for analytic data. No contour extraction, ridge detection, or topological decision is
needed; the only requirement is $\nabla\hat\phi\ne0$ where $\hat\phi$ estimates the function $\phi$
through which $U=G\circ\phi$ factors. If $G$ is not monotone on the band (a ridge), an interpolant of the samples of $U$ has $\nabla\hat U=0$ along the
crest. The crest is then a critical set of its foliation. Thus, a single field-estimate chart cannot cross the crest, and
the band must be cut there (Section~\ref{sec:discussion}). On noisy samples, $\hat\phi$ must be smoothed.
\item \emph{Impose foliation compatibility by construction.} Set $\xi_d=\Lambda(\hat\phi)$ with $\Lambda$ monotone, chosen
by one-dimensional equidistribution across the leaves. Then Theorem~\ref{thm:separability} gives rank one for an exact
factorization $U=G\circ\hat\phi$. Furthermore, Theorem~\ref{thm:leafinjectivity} reduces injectivity of the chart to
injectivity of its leaf charts, with no curvature or reach condition.
\item \emph{Choose the fiber coordinates.} In $d=2$, parametrize each level curve monotonically and carry the
parametrization along the flow of $\hat V$ (Section~\ref{sec:foliation}). Given $\nabla\hat\phi\ne0$ as computed, the
chart is then injective by Theorem~\ref{thm:leafinjectivity}. In $d\ge3$, the leaf charts are surface parametrizations,
subject to the topological obstruction and the open problem discussed in Section~\ref{sec:discussion}.
\item \emph{Account for the error.} The only approximation is how well $U$ factors through $\hat\phi$, measured by the
leafwise defect $\eta$. Sections~\ref{sec:misalign} and~\ref{sec:foliation} bound what it costs.
\end{enumerate}

The construction is one of two options, and the statement below records when each is preferable. It is empirical, and
Section~\ref{sec:numerics} reports the comparison.

\begin{empirical}[Design principle (empirical): when to estimate the field, and when a curve]\label{emp:chartchoice} We consider two cases: 

\emph{Parallel level sets.} Suppose that the level sets of $U$ are parallel, as for a front of constant width or a ridge.
Then an offset chart on a smoothed estimate of one level curve confines the estimation error to one function of one
variable, the normal displacement of that curve; for a straight level curve, the defect then has two separable terms
(Proposition~\ref{note:prop:tilt}). The chart is accurate, of rank comparable to that of a field-estimate chart,
and the cheapest to store. Moreover, at a crest, it is the only option.

\emph{Nonparallel level sets.} Suppose that the level sets are not parallel, as for a front of varying width or the wavy
radial field. Then no offset chart is foliation-compatible, even on the exact curve, while a field-estimate chart is. A
chart built from a field estimate then reaches the accuracy which the data support where offset charts stall, at lower or
equal rank.

The statement holds under four conditions, which are part of it:
\begin{enumerate}[label=(\roman*)]
\item the tolerance is matched to the accuracy of the data;
\item the chart is not reused below that accuracy;
\item the chart needs the field estimate $\hat\phi$ to be built; at noise-limited accuracy, it can then be stored
compactly, but not at clean-data accuracy;
\item $\nabla U$ is bounded away from zero on the band.
\end{enumerate}

The tolerance matched to the data was chosen after the experiments. Under the comparison rule fixed beforehand, with
tolerance $10^{-3}$, the outcome for nonparallel level sets depended on the coverage test of the inverse of the
field-estimate chart. With that test as corrected after the run, the outcome was mixed: the field-estimate chart was better
in four of six cases and comparable in two. With the original test, the outcome was not supported. That rule's tolerance
floor of $10^{-3}$ hid the comparison at the accuracy of the data.
\end{empirical}

\section{Numerical experiments}\label{sec:numerics}\label{note:sec:status}

This section reports the numerical evidence for the statements of Sections~\ref{sec:misalign}--\ref{sec:invert} which are
not proved. There are three experiments. First, we measure how far from sharp the bounds of Section~\ref{sec:misalign}
are. Second, we test whether the dichotomy of Theorem~\ref{thm:rankinjdichotomy}, proved for a straight centerline, holds on
curved ones. Third, we compare
field-estimate and offset charts, which is the evidence for Empirical statement~\ref{emp:chartchoice}.

\subsection{Methods}\label{sec:num:methods}

Every check of a theorem samples the analytic field at the nodes of the chart. Thus, no resampling enters the matrix whose
rank is measured. In the comparison of charts, the data are given on a grid. There, they are evaluated off the grid by a
type-2 nonuniform fast Fourier transform rather than by local interpolation, since the error of local interpolation is
far above the tolerances of the comparison at clean data: $9.7\times10^{-5}$ for linear and $4.7\times10^{-7}$ for cubic interpolation at
resolution $n=801$. %

\emph{Ranks and tolerances.} All ranks are relative Frobenius $\varepsilon$-ranks (Definition~\ref{def:epsrank}), at the
tolerances stated with each result. A rank of at least $0.9\min(N,M)$ is called saturated. It is reported as a lower
bound, drawn with an upward arrow, and, in the rank-cost experiments, recomputed on a larger grid.

\emph{Pre-registration.} The three experiments were pre-registered. That is, their predictions and comparison rules were
fixed and recorded before any run, and every result is reported against them, including the predictions which did not
hold. Where a choice was made after a run, we say so. Two such choices change a reading, both in the comparison of
charts: the noise-matched tolerance, and a correction of the coverage test of the inverse of the field-estimate chart
(Section~\ref{sec:num:chartchoice}). %

\subsection{Rank cost against the rank of the defect}\label{sec:num:rankcost}

Figure~\ref{fig:e4} (Section~\ref{sec:misalign}) shows the $\varepsilon$-rank of $\widetilde U_{ij}=G(\bar\Psi_j+\eta_{ij})$
against the rank $q$ of the defect, for $G(z)=e^{-z^2/w^2}$, $w=0.12$, on the base grid. The defects are rescaled to the
amplitudes $\eta_\infty\in\{3\times10^{-4},3\times10^{-3},10^{-2}\}$, with ten seeds per configuration;
Figure~\ref{fig:e4} shows $\eta_\infty=3\times10^{-3}$, and Table~\ref{tab:e4r} covers all three. They come in four families.
In three of them, the defect is structured in the level variable: it is a polynomial in the level variable with white-noise fiber-side
factors, or a sum of $q$ separable smooth modes with equal or with halving amplitudes. In the fourth, the factors have
independent Gaussian entries.

\emph{Structure, not rank, decides the cost.} For $q\ge3$, the median ranks are ordered polynomial $\le$ smooth with equal
amplitudes $\le$ independent Gaussian factors, at both tolerances. At $\varepsilon=10^{-10}$, $q=20$ and $\eta_\infty=3\times10^{-3}$, the observed
medians are $23$ for the polynomial family and $70.5$ for the smooth family with equal amplitudes. On the other hand, the
family with independent Gaussian factors saturates the base grid at $401$, and it has rank $1023$ of a possible $1024$ on
the large grid. Thus, at equal rank and amplitude, the structure of the defect in the level variable decides the cost, as
Section~\ref{sec:misalign} asserts.

\emph{How sharp the bounds are.} For a random polynomial $G$ of degree $m$ and generic factors, all $60$ cases have
numerical rank exactly $\binom{q+m}{m}$, as Proposition~\ref{prop:binomialattained} states. For the Gaussian profile,
Table~\ref{tab:e4r} gives the ratio of each bound to the observed rank. Over the unsaturated cells of the base grid, the
refined bound is at most $6.7$ times the observed median, while the crude bound is up to $1.2\times10^{5}$ times the
observed median. The refined bound is computed from the measured defect, so it is an a posteriori bound. No bound is
violated in any check in which it is below $\min(N,M)$, as it must be.

\begin{table}[tbhp]
\centering\scriptsize
\caption{Bound over observed rank (median over ten seeds) on the unsaturated cells of the base grid, $q\in\{1,2,3,5,10,20\}$,
$\varepsilon\in\{10^{-6},10^{-10}\}$, $\eta_\infty\in\{3\times10^{-4},3\times10^{-3},10^{-2}\}$: median and maximum per family, for the
refined bound of Theorem~\ref{thm:rankcost}, the bound of Theorem~\ref{thm:productclosed} and the crude bound.}\label{tab:e4r}
\begin{adjustbox}{max width=\textwidth}\begin{tabular}{lrrrrrr}
\toprule
family & cells & refined/obs.\ median & refined/obs.\ max & product-closed bound/obs.\ median & product-closed bound/obs.\ max & crude/obs.\ max \\
\midrule
poly-in-level & 36 & 2.5 & 5.0 & 3.9 & 26.7 & 119\,542 \\
smooth-equal & 36 & 2.5 & 5.4 & 3.7 & 16.0 & 34\,923 \\
smooth-decay & 36 & 2.5 & 5.5 & 4.3 & 44.6 & 97\,128 \\
iid & 33 & 2.2 & 6.7 & -- & -- & 182 \\
\bottomrule
\end{tabular}
\end{adjustbox}
\end{table}

\subsection{The dichotomy beyond straight centerlines}\label{sec:num:dichotomy}

Theorem~\ref{thm:rankinjdichotomy} is proved for a straight centerline. We tested whether the same two behaviors hold
on curved ones, on three centerlines, namely, a straight one, a circle, and a sinusoidal stripe:
$\gamma(t)=(t,0)$ with period $L_{\mathrm{per}}=3$, the circle of radius $R_0=1$, and $\gamma(t)=(t,a\sin\omega t)$ with
$a=0.175$, $\omega=2\pi/3$, $L_{\mathrm{per}}=3$. Each is perturbed
by $\delta=A\cos(2\pi K\theta/L_{\mathrm{per}}+\varphi_0)$ with $A\in\{10^{-3},3\times10^{-3}\}$ and $K\le64$. Every $K$
for which $\mu\in[0.5,2]$ is added, so that $\mu=1$ is bracketed. This gives $156$ configurations on $8192\times401$
nodes, and Figure~\ref{fig:e3rank} (Section~\ref{sec:invert}) shows the rank against $\mu$.

In all $151$ configurations with $|\mu-1|\ge0.02$, the Jacobian is positive and the foot-point map single-valued if and
only if $\mu<1$. On the other hand, the $\varepsilon$-rank of the exact chart at $\varepsilon=10^{-8}$ grows by a factor
of at most $2.25$ from $K=1$ to $K=64$, while $\mu$ grows up to $4096$-fold. For $\mu<1$, it never exceeds the rank $4$ of
the leading-order model by more than $2$. Table~\ref{tab:inject} (Section~\ref{sec:invert}) is a separate run at the
smaller offset $s_{\max}=3w/\sqrt2$, and it shows the same two behaviors.

\subsection{When to estimate the field, and when a curve}\label{sec:num:chartchoice}

\begin{figure}[tbhp]
\centering
\includegraphics{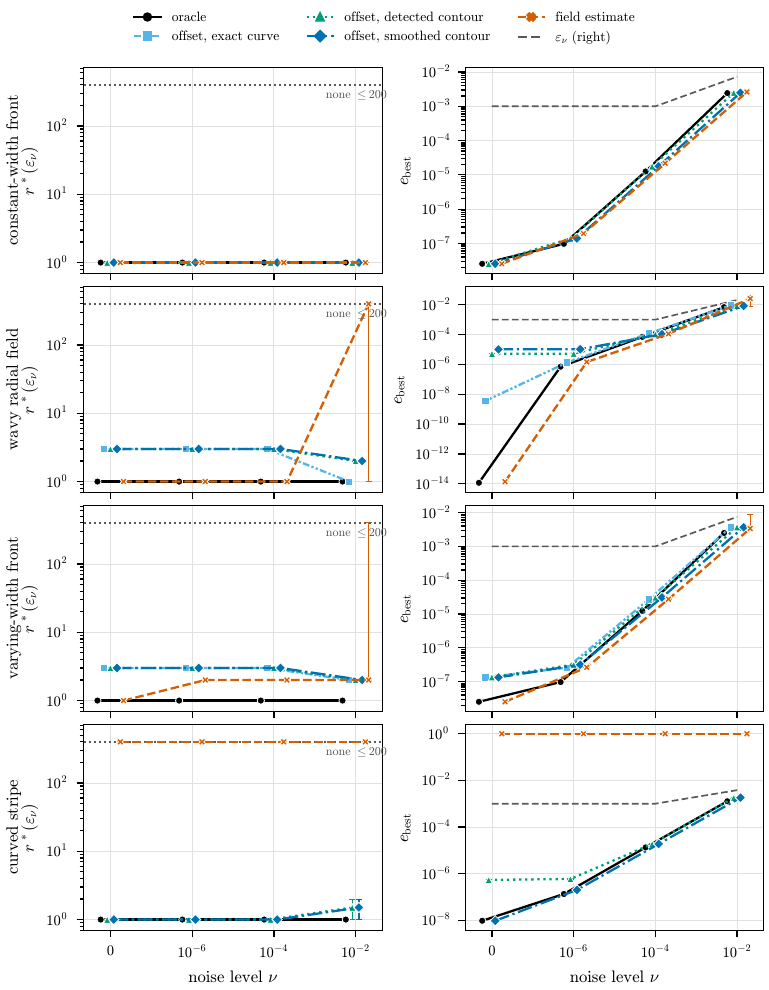}
\caption{$r^*(\varepsilon_\nu)$ (least rank with end-to-end error $\le\varepsilon_\nu$, $r\le200$) and $e_{\rm best}$ against the noise level $\nu$ for the oracle,
the offset chart of the exact curve, offset charts of the detected and of the smoothed contour, and the field-estimate chart, each at its best smoothing width;
median and min/max over seeds (ten for $\nu>0$). $1024\times401$ chart nodes; fields: a constant-width front, the wavy radial field, a varying-width front, and a curved stripe with a crest
(no field-estimate chart). $\varepsilon_\nu=\max(10^{-3},3\min e_{\rm best})$ (the pre-registered rule; its floor hides the noise-free comparison, see Table~\ref{tab:n11tau}).}\label{fig:n1main}
\end{figure}

Empirical statement~\ref{emp:chartchoice} rests on a comparison of charts on four fields, each of which lives on a cylinder
around a sinusoidal curve. These are a front of constant width, which has parallel level sets; a front whose width varies
by $\pm30\%$ along it; a curved Gaussian stripe with a crest; and the wavy radial field $U=e^{-\phi^2}$,
$\phi=|x|+a_w\cos(K\vartheta)$, $a_w=0.25$, $K=5$. The data are $1024^2$ samples with independent Gaussian noise of relative size
$\nu\in\{0,10^{-6},10^{-4},10^{-2}\}$, with ten seeds for $\nu>0$. Every chart is built from the data smoothed by a
Gaussian filter of width $h_f\in\{0,1,2,4,8\}h$, and the data are evaluated at its $1024\times401$ nodes. Five charts are
compared, each at its best smoothing width. These are the oracle, which follows the exact level sets; the offset chart of
the exact central level curve; offset charts of a detected contour and of a smoothed contour; and the field-estimate chart
of Section~\ref{sec:invert}. The statistic is the least rank $r^*$ whose reconstruction error meets the tolerance
(Figure~\ref{fig:n1main}).

\emph{Noise-free data.} The charts which follow the level sets have rank one, as Theorem~\ref{thm:separability} predicts.
For example, the oracle has $\sigma_2/\sigma_1\le3.78\times10^{-14}$. The offset chart of the exact central curve is not
foliation-compatible when the level sets are not parallel. Its $\varepsilon$-rank at $10^{-10}$ is $8$ on the wavy radial
field and $9$ on the front of varying width.

\emph{The pre-registered comparison.} Under the rule fixed before the runs, the tolerance is
$\varepsilon_\nu=\max(10^{-3},3\min e_{\rm best})$, and two ranks are comparable unless they differ by a factor of two.
The rule judges the three fields without a crest at $\nu\le10^{-4}$: the front of constant width tests the parallel case,
and the other two the nonparallel one; the curved stripe and $\nu=10^{-2}$ are not judged (Table~\ref{tab:n1class}). For parallel level sets, every judged cell is comparable. For nonparallel level sets, the field-estimate chart
is better in four of the six judged cells and comparable in two, and no judged cell favors the offset chart. This outcome uses a coverage
test of the inverse which was corrected after the first run. With the original test, three of the six cells favor the
offset chart. Thus, the outcome is mixed with the corrected test, and
not supported with the original one. Moreover, the floor $10^{-3}$ decides every judged cell and hides the comparison at
the accuracy of the data. Therefore, Table~\ref{tab:n11tau} repeats the comparison with the noise-matched tolerance
$\varepsilon_{\rm new}=3e_{\rm best}$, without the floor. This tolerance was chosen after the experiments, and Empirical
statement~\ref{emp:chartchoice} quotes its ranks with that label.

\emph{The mechanism.} At strict tolerances, the picture reverses on noisy data. At $\nu=10^{-4}$ on the front of varying
width, the field-estimate chart has $\varepsilon$-rank $24$ at $10^{-6}$ and $72$ at $10^{-10}$, against $8$ and $16$ for
the smoothed-contour chart. The offset chart of a
smoothed contour carries an error which is a function of the fiber variable. For example, on the front of constant width,
$\sigma_3(\eta)/\sigma_1(\eta)$ is at most $1.4\times10^{-6}$ for $\nu\le10^{-4}$. On the other hand, the noise in
$\hat\phi$ gives the field-estimate chart a defect whose level dependence is rough. In one example
($\nu=10^{-4}$, seed $0$), $\sigma_2(\eta)/\sigma_1(\eta)=0.92$ against $7.0\times10^{-6}$. Thus, the
rank follows the level dependence of the defect, not its amplitude, as Sections~\ref{sec:misalign} and~\ref{sec:foliation} predict.

\emph{The crest and the cost of the chart.} On the curved stripe, a single field-estimate chart fails to reach the crest
leaf in all $155$ configurations, as step~C1 of Section~\ref{sec:invert} anticipates. Offset charts of a detected contour,
however, cover it injectively. Condition~(iv) of Empirical statement~\ref{emp:chartchoice} is needed as well: on a wider
band of the wavy radial field, reaching toward the maximum of $U$, the field-estimate chart loses its advantage.
Finally, the comparison counts ranks, not storage. An offset chart stores a curve, while
the field-estimate chart needs the $1024^2$ array behind $\hat\phi$. We make no compression claim beyond the storage figure
quoted in Section~\ref{sec:discussion}.

\begin{table}[tbhp]
\centering\scriptsize
\caption{The noise-matched tolerance $\varepsilon_{\rm new}=3e_{\rm best}$ (no $10^{-3}$ floor; post hoc) and the resulting $r^*$ and labels, against the pre-registered tolerance $\varepsilon_\nu$ and labels; best $h_f$ unchanged. The rows with $\nu=10^{-2}$ and those of the
curved stripe are not judged under the pre-registered rule.}\label{tab:n11tau}
\begin{adjustbox}{max width=\textwidth}\begin{tabular}{lrrrlrrll}
\toprule
field & $\nu$ & $\varepsilon_\nu$ (pre-reg.) & $\varepsilon_{\rm new}=3e_{\rm best}$ & offset & $r^*$, offset & $r^*$, field estimate & label (pre-reg.) & label new (post hoc) \\
\midrule
constant-width front & 0 & 0.001 & $7.6\times10^{-8}$ & smoothed contour & 1 $\to$ 2 & 1 $\to$ 1 & comparable & field estimate better \\
constant-width front & 1e-06 & 0.001 & $2.9\times10^{-7}$ & smoothed contour & 1 $\to$ 2 & 1 $\to$ 3 & comparable & comparable \\
constant-width front & 0.0001 & 0.001 & $3.8\times10^{-5}$ & smoothed contour & 1 $\to$ 2 & 1 $\to$ 3 & comparable & comparable \\
constant-width front & 0.01 & 0.0073 & $7.3\times10^{-3}$ & smoothed contour & 1 $\to$ 1 & 1 $\to$ 1 & comparable & comparable \\
wavy radial field & 0 & 0.001 & $3.5\times10^{-14}$ & smoothed contour & 3 $\to$ $\infty$ & 1 $\to$ 1 & field estimate better & field estimate better \\
wavy radial field & 1e-06 & 0.001 & $2.1\times10^{-6}$ & smoothed contour & 3 $\to$ $\infty$ & 1 $\to$ 3 & field estimate better & field estimate better \\
wavy radial field & 0.0001 & 0.001 & $2.1\times10^{-4}$ & smoothed contour & 3 $\to$ 3 & 1 $\to$ 1 & field estimate better & field estimate better \\
wavy radial field & 0.01 & 0.021 & $2.1\times10^{-2}$ & smoothed contour & 2 $\to$ 2 & $\infty$ $\to$ $\infty$ & offset better & offset better \\
varying-width front & 0 & 0.001 & $7.6\times10^{-8}$ & smoothed contour & 3 $\to$ $\infty$ & 1 $\to$ 1 & field estimate better & field estimate better \\
varying-width front & 1e-06 & 0.001 & $2.9\times10^{-7}$ & smoothed contour & 3 $\to$ $\infty$ & 2 $\to$ 6 & comparable & field estimate better \\
varying-width front & 0.0001 & 0.001 & $3.7\times10^{-5}$ & smoothed contour & 3 $\to$ 5 & 2 $\to$ 3 & comparable & comparable \\
varying-width front & 0.01 & 0.0075 & $7.5\times10^{-3}$ & smoothed contour & 2 $\to$ 2 & 2 $\to$ 2 & comparable & comparable \\
curved stripe & 0 & 0.001 & $2.9\times10^{-8}$ & smoothed contour & 1 $\to$ 3 & $\infty$ $\to$ $\infty$ & offset better & offset better \\
curved stripe & 1e-06 & 0.001 & $4.0\times10^{-7}$ & smoothed contour & 1 $\to$ 2 & $\infty$ $\to$ $\infty$ & offset better & offset better \\
curved stripe & 0.0001 & 0.001 & $4.0\times10^{-5}$ & smoothed contour & 1 $\to$ 2 & $\infty$ $\to$ $\infty$ & offset better & offset better \\
curved stripe & 0.01 & 0.0038 & $3.8\times10^{-3}$ & smoothed contour & 1.5 $\to$ 1.5 & $\infty$ $\to$ $\infty$ & offset better & offset better \\
\bottomrule
\end{tabular}
\end{adjustbox}
\end{table}

\FloatBarrier

\section{Discussion}\label{sec:discussion}\label{note:sec:cost}\label{note:sec:open}

\paragraph{The cost of the chart}
The representation stores the chart, the factors of the sampled field, and whatever evaluating the inverse requires
(Section~\ref{sec:setting}). On an $n\times n$ grid, a chart stored as two dense nodal arrays costs $2n^2$ numbers, which is
more than the raw field. Thus, a densely stored chart charged to a single field cannot pay for itself. A chart shared by $N_f$
fields with the same ranks pays for itself against Cartesian low-rank storage exactly when $N_f\,\Delta r>n$, where $\Delta r$
is the rank it saves per field. Therefore, a chart must be stored compactly, or shared. In our experiments, the first was
possible at noise-limited accuracy: a truncated singular value decomposition compressed the field-estimate chart to about
$1.3$--$1.6\times10^{4}$ numbers, though not at clean-data accuracy (Table~\ref{tab:n12}). Choosing $X$ is
then a design problem, namely, to minimize $C_{\rm tensor}(U\circ X)+\lambda\,C_{\rm map}(X)$ at a prescribed reconstruction
error. This paper bounds the first term only.

\paragraph{Limits of a single chart}
The flow chart of Theorem~\ref{thm:bandproduct} needs a band without critical points of $\phi$. Near a critical point it
breaks down, and offset charts suffer too: near a nondegenerate critical point the level sets curve sharply, so their reach,
and with it the admissible offset, shrinks to zero.
A crest of $U$ is a critical set of another kind, where $G'=0$ while $\nabla\phi\ne0$. A single field-estimate chart cannot
cross it, because $\nabla\hat U=0$ there (step~C1 of Section~\ref{sec:invert}), but an offset chart can when the level sets of $\phi$ are
parallel, since its leaves then are these level sets (Section~\ref{sec:num:chartchoice}). In general, the domain must therefore be cut into
regular bands at the critical values of $\phi$, which contour trees and Reeb graphs compute~\cite{Carr2003,Tierny2018}. In
dimension $d\ge3$, there is a second, topological limit. A smooth, nonsingular tensor-product parametrization of a leaf would
give $d-1$ commuting, linearly independent vector fields on it, and a closed connected manifold which carries such fields
is a torus~\cite{Arnold1989}. In particular, the spheres around a radial $\phi$ admit none; for $d=3$, this already follows
from the hairy ball theorem~\cite{Milnor1965}. For $d=2$, the
leaves are circles, and there is no such obstruction. Singular parametrizations always exist and still give rank one, but
they make the inverse $X^{-1}$ multivalued on their singular set.

\paragraph{Open problems}
(i) \emph{The analyticity radius.} Corollary~\ref{cor:georank} reduces the rank cost to the radius $R_\phi$. However,
$R_\phi$ is known in closed form only for the straight centerline (order two) and for the circle, where it is the distance to
the center of curvature, the critical point of $\phi$. A lower bound on $R_\phi$ from
curvature and reach data, valid for general curved foliations, is missing. The circle suggests the distance from the
estimated leaf to the nearest critical or focal point of $\phi$.

(ii) \emph{Lower bounds and sharpness.} Proposition~\ref{prop:binomialattained} shows that the binomial count is attained
for polynomial profiles. However, there is no lower bound on $r_\varepsilon$ in terms of $R_\phi$, and none for the least
$\varepsilon$-rank over a constrained class of charts (nonsingular, of bounded distortion, or of bounded storage). Without a
constraint, that least rank is one on every regular band (Proposition~\ref{prop:rankone}). In the other direction, the
refined bound of Theorem~\ref{thm:rankcost} is at most $6.7$ times the observed $\varepsilon$-rank in
Section~\ref{sec:numerics}. Quantifying the spectral decay of the Hadamard powers of a matrix of given rank and spectrum would
close that gap.

(iii) \emph{Leaf charts in dimension three and more.} Theorem~\ref{thm:leafinjectivity} reduces injectivity to the leaf
charts. In $d\ge3$, these are surface parametrizations, subject to the topological obstruction described above. A moving mesh on each leaf, with alignment taken with respect to the metric induced
on the leaf, is where variational mesh adaptation applies. Constructing leaf charts which are nonsingular and depend
smoothly on the level is open.

(iv) \emph{Atlases and evolving fields.} Beyond a single band, one needs an atlas of bands and, for $d\ge3$, fiber patches. How rank is accounted across an atlas, and how bands which meet the boundary are treated, is
open. For evolving fields, the comparison with the coordinate flows of~\cite{DektorVenturi2023,DektorVenturi2025} needs
test cases which a linear flow cannot align, such as a front of variable curvature or a front which deforms under shear.

\paragraph{Limitations}
The bounds of Sections~\ref{sec:misalign} and~\ref{sec:foliation} concern the (fiber)\,$\times$\,(level) unfolding of the
sampled field, which for $d\ge3$ is one tensor-train unfolding. The other tensor-train ranks and the CP rank are treated only
in the exact case (Section~\ref{sec:setting}). The numerical experiments are in dimension two, and most use Gaussian profiles.
The noise-matched tolerance of the comparison of charts, and a correction of its coverage test, were chosen after the
experiments; under the pre-registered rule, the outcome for nonparallel level sets is mixed with the corrected test and not
supported with the original one (Section~\ref{sec:numerics}). %
Finally, the theory is for a single regular band, with the limits described above.

\section*{Acknowledgments}
The authors gratefully acknowledge the support of the Laboratory Directed Research and Development (LDRD) program of Los
Alamos National Laboratory under project number 20260495ER. Los Alamos National Laboratory is operated by Triad National
Security, LLC, for the National Nuclear Security Administration of U.S. Department of Energy (Contract
No.~89233218CNA000001).

\paragraph{Use of AI tools}
AI tools were used to draft text, to write and run the code for the numerical experiments, and to check
citations. The authors read every cited work, checked every proof, and verified every number in the paper against the code package.
The authors assume responsibility for all content.

\appendix
\section{Proofs}\label{app:proofs}

In this appendix, we give the proofs of Lemma~\ref{lem:obliqueridge}, Theorem~\ref{thm:rankcost}, and
Corollary~\ref{cor:levelanalytic}.

\subsection{Proof of Lemma~\ref{lem:obliqueridge}}\label{app:proof:lem:obliqueridge}

\begin{proof}[Proof of Lemma~\ref{lem:obliqueridge}]
The idea is that a kernel of finite rank forces $G$ to satisfy a linear ordinary differential equation (ODE) with
constant coefficients. We argue by contradiction. Suppose that the integral operator $T_f:L^{2}(I_2)\to L^{2}(I_1)$ has
finite rank, with range $V$.

\emph{Step 1: the derivatives of the columns lie in $V$.} Apply $T_f$ to a sequence of bumps converging to the Dirac mass
at $t_0$. Since $f$ is continuous and $V$ is closed, this gives $f(\cdot,t_0)\in V$. Furthermore, $t\mapsto f(\cdot,t)$ is
smooth with values in $V$. Thus, its $t$-derivatives also take values in $V$, being limits of difference quotients. Hence,
$s\mapsto b^{j}G^{(j)}(as+bt_0)$ lies in $V$ for every $j$.

\emph{Step 2: $G$ is an exponential polynomial.} Any $\dim V+1$ of these functions are linearly dependent. Thus, $G$
satisfies a nontrivial constant-coefficient linear ODE on the nondegenerate interval $aI_1+bt_0$. By analyticity, it
satisfies the same ODE on $\RR$, hence it is an exponential polynomial. This contradicts the hypothesis on $G$.

\emph{Step 3: the $\varepsilon$-rank tends to infinity.} Suppose that $r_\varepsilon(f)\le r$ for all $\varepsilon>0$.
Then $\sigma_{r+1}(T_f)\le\varepsilon\|f\|$ for every $\varepsilon$. That is, $T_f$ has rank $\le r$, a contradiction.
\end{proof}

\subsection{Proof of Theorem~\ref{thm:rankcost}}\label{app:proof:thm:rankcost}

\begin{proof}[Proof of Theorem~\ref{thm:rankcost}]
The idea is to truncate a Chebyshev expansion in the defect, and then to count the ranks of the Hadamard powers of $\eta$
which the truncation produces.

\emph{Step 1: the Chebyshev truncation.} Fix $j$, and set $f_j(u)=G(\bar\Psi_j+u)$. Consider the (open, filled) Bernstein
ellipse
\[
E_\rho=\{\eta_\infty(w+w^{-1})/2:\ 1/\rho<|w|<\rho\}.
\]
Its semi-major axis is $\eta_\infty(\rho+\rho^{-1})/2=R$ for the $\rho$ of \eqref{note:eq:rho-m}. Thus, it lies in the disk
$|u|<R$, and $f_j$ is analytic in $E_\rho$ with $|f_j|\le M_G$.

We write $f_j=\sum_kc_k(j)T_k(u/\eta_\infty)$. By the theorem on Chebyshev coefficients of analytic
functions~\cite[Thm.~8.1]{Trefethen2019}, we have $|c_k(j)|\le2M_G\rho^{-k}$. Since $|T_k|\le1$ on $[-1,1]$, it follows that
\[
\sup_{|u|\le\eta_\infty}\Bigl|f_j(u)-\sum_{k\le m}c_k(j)T_k(u/\eta_\infty)\Bigr|
\le\sum_{k>m}2M_G\rho^{-k}=\frac{2M_G\rho^{-m}}{\rho-1}\le\frac{\varepsilon'}2 .
\]
Hence, the matrix
\[
\widetilde U_m:=\sum_{k=0}^{m}(\mathbf 1c_k^{\mathsf T})\circ T_k(\eta/\eta_\infty)
\]
satisfies $\|\widetilde U-\widetilde U_m\|_F\le\sqrt{NM}\,\varepsilon'/2=\varepsilon\|\widetilde U\|_F/2$.

\emph{Step 2: the crude bound.} Since $T_k$ is a polynomial of degree $k$, we have
$T_k(\eta/\eta_\infty)=\sum_{p\le k}t_{kp}\eta_\infty^{-p}\,\eta^{\circ p}$. Collecting powers gives
\[
\widetilde U_m=\sum_{p=0}^{m}(\mathbf 1w_p^{\mathsf T})\circ\eta^{\circ p},\qquad
w_p=\sum_{k\ge p}t_{kp}\eta_\infty^{-p}c_k .
\]
Since $(\mathbf 1w^{\mathsf T})\circ B=B\operatorname{diag}(w)$, the $p$th summand has rank at most
$\operatorname{rank}\eta^{\circ p}$. By Lemma~\ref{note:lem:hadamard}, this rank is at most $\binom{q+p-1}{p}$. Furthermore,
$\sum_{p=0}^{m}\binom{q+p-1}{p}=\binom{q+m}{m}$ by the hockey-stick identity. Finally, by Step~1,
$r_\varepsilon(\widetilde U)\le\operatorname{rank}\widetilde U_m$.

\emph{Step 3: the refined bound.} We replace each $T_k(\eta/\eta_\infty)$ by an optimal rank-$q_k$ approximation $B_k$.
Since $\|(\mathbf 1c_k^{\mathsf T})\circ E\|_F\le\|c_k\|_\infty\|E\|_F$, the added Frobenius error is at most
$\sum_k\|c_k\|_\infty\upsilon_k=\varepsilon\|\widetilde U\|_F/2$. Together with Step~1, the total error is at most
$\varepsilon\|\widetilde U\|_F$. Moreover, $\operatorname{rank}\sum_k(\mathbf 1c_k^{\mathsf T})\circ B_k\le\sum_kq_k$.
\end{proof}

\subsection{Proof of Corollary~\ref{cor:levelanalytic}}\label{app:proof:cor:levelanalytic}

\begin{proof}[Proof of Corollary~\ref{cor:levelanalytic}]
The idea is to replace the defect by its truncated Chebyshev series in the level variable, and then to apply
Theorem~\ref{note:thm:structured}(a).

We truncate each $\eta_i$ at its Chebyshev polynomial of degree $n_\upsilon-1$. By the estimate
of~\cite[Thm.~8.1]{Trefethen2019}, summed over $k\ge n_\upsilon$, the entrywise error is at most $\upsilon$. Thus,
$\eta_\upsilon$ is a polynomial of degree $<n_\upsilon$ in $x_j$ with $\max|\eta-\eta_\upsilon|\le\upsilon$.

Entrywise, the points $\bar\Psi_j+\eta_{ij}$ and $\bar\Psi_j+\eta_{\upsilon,ij}$ lie in the interval on which $|G'|\le M_1$.
Thus, by the mean value theorem,
\[
\|G(\bar\Psi+\eta)-G(\bar\Psi+\eta_\upsilon)\|_F\le M_1\upsilon\sqrt{NM}\le\varepsilon\|\widetilde U\|_F/4.
\]

Next, Theorem~\ref{note:thm:structured}(a), applied to $\eta_\upsilon$ at tolerance $\varepsilon/4$, gives a matrix of rank
at most the stated bound. This matrix lies within $\tfrac\varepsilon4(1+\tfrac\varepsilon4)\|\widetilde U\|_F$ of
$G(\bar\Psi+\eta_\upsilon)$. The two errors add to less than $\varepsilon\|\widetilde U\|_F$. Finally, both $m$ and
$n_\upsilon$ are $O(\log(1/\varepsilon))$.
\end{proof}

\section*{Supplementary material}
Appendices~\ref{sm:first} onward contain the data and the derivations to which the paper points.

\section{Metric criteria}\label{sm:sec:metric}\label{sm:first}

\subsection{The counterexample numerically}\label{sm:sec:countertable}

Table~\ref{tab:counter} gives the numbers behind the statement after Theorem~\ref{thm:counterexample}.

\begin{table}[tbhp]
\centering\scriptsize
\caption{Theorem~\ref{thm:counterexample} numerically: $\mathcal M_\phi=\operatorname{diag}(1,9)$, $\chi=1$, the linear chart
$X(\xi)=\sqrt\chi\,\mathcal M^{-1/2}Q_\psi\xi$ on a $400\times400$ grid over $\xi\in[0,1]^2$ (endpoints included), $U=\exp(-x_2^2/w^2)$, $w=0.1$.
The residuals are $\max_{\rm nodes}\|J^{\mathsf T}\mathcal MJ-\chi I\|_{\max}$ (alignment, Proposition~\ref{prop:alignconformal} form) and
$\max_{\rm nodes}|\sqrt{\det\mathcal M}\det J-\chi|$, and they are rounding noise. The condition number is $\mathrm{cond}(DX)=3$ for every $\psi$.
Centered: ridge through $X(\tfrac12,\tfrac12)$. At $\varepsilon=10^{-14}$, the ranks are limited by the SVD rounding floor (an exactly
rank-one $400^2$ matrix has computed tail $2.7\times10^{-14}$) and are not shown.}\label{tab:counter}\label{note:tab:counter}
\begin{adjustbox}{max width=\textwidth}\begin{tabular}{rrrrrrrr}
\toprule
$\psi$ (deg) & alignment res. & equidistr. res. & $\mathrm{cond}(DX)$ & $r_{10^{-6}}$ & $r_{10^{-10}}$ & $r_{10^{-12}}$ & centred $r_{10^{-6}}/r_{10^{-10}}$ \\
\midrule
0 & 0 & 0 & 3.000 & 1 & 1 & 1 & 1/1 \\
5 & $1.1\times10^{-16}$ & 0 & 3.000 & 4 & 7 & 8 & 5/7 \\
15 & $2.2\times10^{-16}$ & 0 & 3.000 & 6 & 9 & 10 & 7/10 \\
30 & $1.3\times10^{-16}$ & 0 & 3.000 & 6 & 10 & 12 & 8/11 \\
45 & $1.1\times10^{-16}$ & 0 & 3.000 & 7 & 10 & 12 & 8/12 \\
60 & $4.7\times10^{-17}$ & $2.2\times10^{-16}$ & 3.000 & 6 & 10 & 12 & 8/11 \\
90 & $7.5\times10^{-33}$ & 0 & 3.000 & 1 & 1 & 1 & 1/1 \\
\bottomrule
\end{tabular}
\end{adjustbox}
\end{table}

\subsection{The rotated Gaussian}\label{sm:sec:mehler}

\begin{proof}[Proof of Proposition~\ref{prop:rotgauss}]
Put $a=x_1/\sqrt{\Sigma_{11}}$ and $b=x_2/\sqrt{\Sigma_{22}}$. This diagonal rescaling multiplies all singular values by
one constant. Then $U=\exp\bigl(-(a^{2}-2\varrho ab+b^{2})/(2(1-\varrho^{2}))\bigr)$. Let
$h_n(x)=(2^{n}n!\sqrt\pi)^{-1/2}H_n(x)e^{-x^{2}/2}$ be the Hermite functions, an orthonormal basis of $L^{2}(\RR)$. In
this form, Mehler's formula~\cite{Mehler1866,Watson1933} reads, for $|z|<1$,
\[
\sum_{n\ge0}z^{n}h_n(x)h_n(y)=\frac{1}{\sqrt{\pi(1-z^{2})}}
\exp\Bigl(-\frac{(1+z^{2})(x^{2}+y^{2})-4zxy}{2(1-z^{2})}\Bigr).
\]

It says that the kernel
\[
K_z(x,y)=\exp\bigl(-[(1+z^{2})(x^{2}+y^{2})-4zxy]/(2(1-z^{2}))\bigr)
\]
has the singular value decomposition $K_z=\sqrt{\pi(1-z^{2})}\sum_nz^{n}\,h_n\otimes h_n$, with singular values
$\sqrt{\pi(1-z^{2})}\,|z|^{n}$. Choose $z\in(-1,1)$ with $2z/(1+z^{2})=\varrho$, that is,
$z=\varrho/(1+\sqrt{1-\varrho^{2}})$, and choose $\lambda>0$ with $\lambda^{2}/(1-\varrho^{2})=(1+z^{2})/(1-z^{2})$.
Comparing the coefficients of $x^{2}+y^{2}$ and of $xy$, the substitution $a=\lambda x$, $b=\lambda y$ turns $U$ into
$K_z$, and it multiplies all singular values by the constant $\lambda$. Hence, $\sigma_{k+1}/\sigma_k=|z|$, the ratio of
the statement.

The tail energy is $\sum_{n\ge r}|z|^{2n}\propto |z|^{2r}$. Thus, the relative Frobenius error of the rank-$r$ truncation
is exactly $|z|^{r}$, and $r_\varepsilon=\lceil\ln(1/\varepsilon)/\ln(1/|z|)\rceil$.

For the $45^\circ$ rotation of an axis-aligned Gaussian with variances $\varsigma_1^{2},\varsigma_2^{2}$, one has
$\Sigma_{11}=\Sigma_{22}=(\varsigma_1^{2}+\varsigma_2^{2})/2$ and $\Sigma_{12}=(\varsigma_1^{2}-\varsigma_2^{2})/2$. Since
$\kappa_\Sigma=\varsigma_1^{2}/\varsigma_2^{2}$, this gives
\[
\varrho=\frac{\kappa_\Sigma-1}{\kappa_\Sigma+1},\qquad
\sqrt{1-\varrho^{2}}=\frac{2\sqrt{\kappa_\Sigma}}{\kappa_\Sigma+1},\qquad
|z|=\frac{\sqrt{\kappa_\Sigma}-1}{\sqrt{\kappa_\Sigma}+1}.
\]
Finally, $\ln\bigl((\sqrt{\kappa_\Sigma}+1)/(\sqrt{\kappa_\Sigma}-1)\bigr)\sim2/\sqrt{\kappa_\Sigma}$.
\end{proof}

\subsection{The condition number numerically}\label{sm:sec:whitening}

Table~\ref{tab:whitening} gives the numbers behind Figure~\ref{fig:intro} and Section~\ref{sec:condnumber}.

\begin{table}[tbhp]
\centering\scriptsize
\caption{The condition number of the chart does not see the rank: rotated Gaussian with $\kappa_\Sigma=25$ at $30^\circ$,
sampled on $201\times201$ nodes on the stated $\xi$-boxes, outside which the relative $L^2$ mass of the field is at most
$10^{-10}$.}\label{tab:whitening}
\begin{adjustbox}{max width=\textwidth}\begin{tabular}{llrrrr}
\toprule
chart & $\xi$-box & $\mathrm{cond}(DX)$ & $r_{10^{-6}}$ & $r_{10^{-10}}$ & $\sigma_2/\sigma_1$ \\
\midrule
cartesian & [-4.6,4.6]$\times$[-2.8,2.8] & 1 & 30 & 50 & $6.3\times10^{-1}$ \\
rotation & [-5.2,5.2]$\times$[-1.1,1.1] & 1 & 1 & 1 & $7.0\times10^{-16}$ \\
whitening & [-5.2,5.2]$\times$[-5.2,5.2] & 5 & 1 & 1 & $9.3\times10^{-16}$ \\
\bottomrule
\end{tabular}
\end{adjustbox}
\end{table}

\section{The rank cost: constants for a Gaussian profile}\label{sm:sec:misalign}\label{sm:sec:workedinstance}

The Gaussian profile $G(z)=e^{-z^2/w^2}$ is entire, so any $R$ may be used in Theorem~\ref{thm:rankcost}. With $R=w$, one
may take $M_G=e$, since for $|y|\le w$
\[
|G(x+iy)|=e^{(y^2-x^2)/w^2}\le e^{y^2/w^2}\le e .
\]
These are the constants with which~\eqref{note:eq:rho-m} gives the values of $m$ that Section~\ref{sec:misalign} quotes
for the setting of Figure~\ref{fig:e4}.

\section{The foliation order: two examples}\label{sm:sec:foliation}

\subsection{A gauge which is not flow-transported}\label{sm:sec:gaugeexample}

For example, a defect $\eta=A\sin Ky$ is of rank one and order one in the flow-transported gauge. Consider the gauge
$y=\xi'+b\,\xi_d\sin\xi'$, with $\xi_d$ here the level coordinate normalized to $[-1,1]$. Then the Jacobi--Anger expansion
gives
\[
\eta(\xi_d,\xi')=A\sum_{n\in\mathbb Z}J_n(Kb\,\xi_d)\sin\bigl((K+n)\xi'\bigr).
\]

Furthermore, $|J_n(x)|\le(|x|/2)^{|n|}/|n|!$ for real $x$~\cite[(10.14.4)]{DLMF}. Thus, at most
$O\bigl(Kb\max|\xi_d|+\log(1/\varepsilon)\bigr)$ modes of the level dependence are significant. Here the identity
$J_{-n}=(-1)^nJ_n$ pairs the modes $K\pm n$ onto one level function. That is, the strength of the gauge multiplies the fiber
frequency.

Numerically, the $\varepsilon$-rank of $\widetilde U$ is $4$ for $b=0$ and $40$, $71$, $120$ for $b=0.25$, $0.5$, $0.9$
(for $K=16$: $4$, $16$, $24$, $37$). Here $G(z)=e^{-z^2/w^2}$ with $w=0.12$, on $401$ levels in $[-s_{\max},s_{\max}]$,
$s_{\max}=3w/\sqrt2$, and $4096$ periodic fiber points; the gauge acts on the normalized level coordinate,
$A=3\times10^{-3}$, $K=64$, and $\varepsilon=10^{-8}$. %

\subsection{The circle}\label{sm:sec:circle}

\begin{remark}[The circle]\label{note:rem:circle}\label{rem:circleexample}
Let $\gamma$ be the circle of radius $R_0$, let $\phi=|x|-R_0$, and let $\hat\gamma=(R_0+\delta)e_r$. The offset point is
\[
P=(R_0+\delta)(1+s/\hat\ell)e_r-(s\delta'/\hat\ell)e_\vartheta
\]
with $\hat\ell^2=(R_0+\delta)^2+\delta'^2$. Thus, $\eta=|P|-R_0-s$ with
\[
|P|=\hat\ell^{-1}\sqrt{(R_0+\delta)^{2}(\hat\ell+s)^{2}+s^{2}\delta'^{2}}.
\]

As a function of $s$, this has branch points at $s=-\hat\ell(R_0+\delta)/\bigl((R_0+\delta)\pm i\delta'\bigr)$, that is,
on the circle $|s|=R_0+\delta(\theta)$, \emph{independently of the wavenumber of $\delta$}. Hence, $\eta(\theta,\cdot)$ is
analytic in the disk $|s|<R_0+\delta(\theta)$ about the central leaf $s=0$. The every-leaf bound of
Corollary~\ref{note:cor:geometric} would use the smaller radius $R_0+\min\delta-s_{\max}$, the distance from the edges of
the band to the singular circle. The wavenumber enters only through $M_\eta$, which grows with the slope. As
Section~\ref{sec:invert} shows (Table~\ref{tab:inject}), the exact chart can fold while the rank moves only from $4$ to
$11$.

Numerically, the Chebyshev coefficients of $\eta(\theta,\cdot)$ on $[-s_{\max},s_{\max}]$ decay at the rate $0.10$--$0.13$
for $K=1,16,64,256$. This compares with $\rho_s^{-1}=0.129$, predicted from the central-leaf radius $R_\phi=R_0+\delta$ by a
Bernstein-ellipse count, which needs analyticity of $\eta(\theta,\cdot)$ in the ellipse only. For small slopes, the two
conjugate branch points nearly coincide, and the observed ratio equals $\rho_s^{-1}$. For large slopes, it is lowered by
the oscillation of the conjugate pair as well as by the prefactor $k^{-3/2}$.
\end{remark}

\section{Numerical experiments: the comparison of charts}\label{sm:sec:numerics}\label{sm:sec:chartdata}

Table~\ref{tab:n1class} gives the comparison of Section~\ref{sec:num:chartchoice} under the rule fixed before the runs. For the front of constant width
(parallel level sets), every judged cell is comparable, with $r^*=1$ for both the smoothed-contour offset chart and the
field-estimate chart. For nonparallel level sets, on the wavy radial field, $r^*=1$ against $3$ at all three judged noise
levels; on the front of varying width, $1$ against $3$ on noise-free data and $2$ against $3$ with noise. One definition
was corrected after the first run: the inverse of the field-estimate chart declared a point uncovered when its seed curve
missed the leaf by more than $10^{-6}$, which the jittered seeds of the wavy radial field exceed. With that original test,
$47.0\%$ of its points are falsely uncovered, and three of the six cells favor the offset chart.

\begin{table}[tbhp]
\centering\scriptsize
\caption{Classification under the pre-registered rule (factor 2 on the median $r^*$; revised coverage threshold of the field-estimate chart's inverse, post hoc; the original threshold alongside). Case: the paragraph of Empirical statement~\ref{emp:chartchoice} that the cell tests (parallel or non-parallel level sets), or ``not judged'' for the cells outside the pre-registered comparison.}\label{tab:n1class}
\begin{adjustbox}{max width=\textwidth}\begin{tabular}{llrrlrrll}
\toprule
case & field & $\nu$ & $\varepsilon_\nu$ & offset & $r^*$, offset ($h_f$) & $r^*$, field estimate ($h_f$) & label (revised) & label (original) \\
\midrule
parallel & constant-width front & 0 & 0.001 & smoothed contour & 1 (0) & 1 (0) & comparable & comparable \\
parallel & constant-width front & 1e-06 & 0.001 & smoothed contour & 1 (0) & 1 (4) & comparable & comparable \\
parallel & constant-width front & 0.0001 & 0.001 & smoothed contour & 1 (4) & 1 (8) & comparable & comparable \\
not judged & constant-width front & 0.01 & 0.0073 & smoothed contour & 1 (4) & 1 (8) & comparable & comparable \\
non-parallel & wavy radial field & 0 & 0.001 & smoothed contour & 3 (1) & 1 (0) & field estimate better & offset better \\
non-parallel & wavy radial field & 1e-06 & 0.001 & smoothed contour & 3 (2) & 1 (1) & field estimate better & offset better \\
non-parallel & wavy radial field & 0.0001 & 0.001 & smoothed contour & 3 (8) & 1 (1) & field estimate better & offset better \\
not judged & wavy radial field & 0.01 & 0.021 & smoothed contour & 2 (8) & $\infty$ (2) & offset better & offset better \\
non-parallel & varying-width front & 0 & 0.001 & smoothed contour & 3 (0) & 1 (0) & field estimate better & field estimate better \\
non-parallel & varying-width front & 1e-06 & 0.001 & smoothed contour & 3 (0) & 2 (4) & comparable & comparable \\
non-parallel & varying-width front & 0.0001 & 0.001 & smoothed contour & 3 (0) & 2 (4) & comparable & comparable \\
not judged & varying-width front & 0.01 & 0.0075 & smoothed contour & 2 (4) & 2 (8) & comparable & comparable \\
not judged & curved stripe & 0 & 0.001 & smoothed contour & 1 (0) & $\infty$ (0) & offset better & offset better \\
not judged & curved stripe & 1e-06 & 0.001 & smoothed contour & 1 (1) & $\infty$ (0) & offset better & offset better \\
not judged & curved stripe & 0.0001 & 0.001 & smoothed contour & 1 (4) & $\infty$ (0) & offset better & offset better \\
not judged & curved stripe & 0.01 & 0.0038 & smoothed contour & 1.5 (8) & $\infty$ (0) & offset better & offset better \\
\midrule
\multicolumn{9}{l}{non-parallel: mixed (original threshold: not supported); parallel: supported (original threshold: supported)} \\
\bottomrule
\end{tabular}
\end{adjustbox}
\end{table}

Table~\ref{tab:n12} gives the storage of the field-estimate chart which Section~\ref{sec:discussion} quotes. The
coordinates of the chart are compressed by a truncated singular value decomposition, and the chart is inverted by Newton's
method on the compressed chart. At $\nu=10^{-4}$, this brings the field-estimate chart to about $1.3$--$1.6\times10^{4}$
numbers; at clean-data accuracy, it does not compress.

\begin{table}[tbhp]
\centering\scriptsize
\caption{Storage of the charts with compressed coordinates. Coordinates $x_1,x_2$ of each chart truncated by SVD (position error at most $f$ times the median node spacing,
$f\in\{0.5,\dots,10^{-12}\}$); smallest total $T^*=$ chart $+\,1426\,r^*(\varepsilon_{\rm new})$ with $e_{\rm best}\le1.1\times$ the uncompressed reference. ``self'': inverse by Newton on the
compressed chart (self-contained); ``same$^\dagger$'': the inverse of Figure~\ref{fig:n1main}, which for the field-estimate chart still evaluates $\hat\phi$ (not counted). ``uncompressed total'': the smallest total with the uncompressed chart (for the field-estimate chart, the $1024^2$ array behind $\hat\phi$ and the seeds). Medians over seeds 0--2 at $\nu=10^{-4}$.}\label{tab:n12}
\begin{adjustbox}{max width=\textwidth}\begin{tabular}{lrllrlrrrrr}
\toprule
field & $\nu$ & chart & inverse & ref.\ $e_{\rm best}$ & ranks ($x_1,x_2$) & chart & $e_{\rm best}$ ratio & $r^*$ & total $T^*$ & uncompressed total \\
\midrule
wavy radial field & 0 & field estimate & self & $5.4\times10^{-6}$ & 4,4 & 11\,408 & 1.00 & $\infty$ & $\infty$ & 1\,052\,050 \\
wavy radial field & 0 & field estimate & same$^\dagger$ & $1.4\times10^{-14}$ & 17,17 & 48\,484 & 0.75 & 1 & 49\,910 & 1\,052\,050 \\
wavy radial field & 0 & smoothed contour & self & $5.9\times10^{-8}$ & 2,2 & 5\,704 & 1.00 & $\infty$ & $\infty$ & $\infty$ \\
wavy radial field & 0 & smoothed contour & same$^\dagger$ & $1.0\times10^{-5}$ & 2,2 & 5\,704 & 1.00 & $\infty$ & $\infty$ & $\infty$ \\
wavy radial field & 0.0001 & field estimate & self & $1.1\times10^{-4}$ & 4,4 & 11\,408 & 0.99 & 1 & 12\,834 & 1\,052\,050 \\
wavy radial field & 0.0001 & field estimate & same$^\dagger$ & $1.1\times10^{-4}$ & 4,4 & 11\,408 & 1.00 & 1 & 12\,834 & 1\,052\,050 \\
wavy radial field & 0.0001 & smoothed contour & self & $1.0\times10^{-4}$ & 2,2 & 5\,704 & 1.00 & 3 & 9\,982 & 7\,356 \\
wavy radial field & 0.0001 & smoothed contour & same$^\dagger$ & $1.1\times10^{-4}$ & 2,2 & 5\,704 & 1.00 & 3 & 9\,982 & 7\,356 \\
varying-width front & 0 & field estimate & self & $2.2\times10^{-7}$ & 7,7 & 19\,964 & 1.00 & $\infty$ & $\infty$ & 1\,052\,050 \\
varying-width front & 0 & field estimate & same$^\dagger$ & $2.5\times10^{-8}$ & 64,72 & 193\,936 & 1.08 & 1 & 195\,362 & 1\,052\,050 \\
varying-width front & 0 & smoothed contour & self & $1.3\times10^{-7}$ & 2,2 & 5\,704 & 1.00 & $\infty$ & $\infty$ & $\infty$ \\
varying-width front & 0 & smoothed contour & same$^\dagger$ & $1.3\times10^{-7}$ & 2,2 & 5\,704 & 1.00 & $\infty$ & $\infty$ & $\infty$ \\
varying-width front & 0.0001 & field estimate & self & $2.7\times10^{-5}$ & 3,4; 4,4 & 9\,982 & 1.05 & 4 & 15\,686 & 1\,054\,902 \\
varying-width front & 0.0001 & field estimate & same$^\dagger$ & $2.7\times10^{-5}$ & 19,19; 20,21; 21,22 & 58\,466 & 1.01 & 3 & 62\,744 & 1\,054\,902 \\
varying-width front & 0.0001 & smoothed contour & self & $3.1\times10^{-5}$ & 2,2 & 5\,704 & 1.00 & 5 & 12\,834 & 10\,208 \\
varying-width front & 0.0001 & smoothed contour & same$^\dagger$ & $3.1\times10^{-5}$ & 2,2 & 5\,704 & 1.00 & 5 & 12\,834 & 10\,208 \\
\bottomrule
\end{tabular}
\end{adjustbox}
\end{table}

\section{Discussion: reach at a nondegenerate critical point}\label{sm:sec:limits}\label{sm:sec:reachproof}

The paragraph \emph{Limits of a single chart} of Section~\ref{sec:discussion} uses the following proposition.

\begin{proposition}[Reach degenerates at nondegenerate critical points]\label{note:prop:reach-crit}\label{prop:reachcritical}
Let $\phi\in C^{2}(\Omega)$, $d\ge2$, have a nondegenerate critical point $x_0$ with value $c$. Let $\alpha_{\min}>0$ be the
smallest absolute value of an eigenvalue of $H=D^{2}\phi(x_0)$. Consider each sign of ${t}$ for which $H$ has an
eigenvalue of that sign (one sign only at a local extremum), and let $|{t}|$ be small. Then the level set
$\Sigma_{t}=\phi^{-1}(c+{t})$ contains a point $x_{t}$ at which the second fundamental form
$\II=-P_nD^{2}\phi\,P_n/|\nabla\phi|$ of $\Sigma_{t}$ satisfies $\|\II(x_{t})\|\ge c_0|{t}|^{-1/2}$. Here $P_n$ is the
orthogonal projection onto the tangent space of $\Sigma_{t}$, and $c_0>0$ does not depend on ${t}$. Consequently, if $c+{t}$
is a regular value and $\Sigma_{t}$ is a compact hypersurface without boundary, then $\rch(\Sigma_{t})=O(\sqrt{|{t}|})$.
\end{proposition}

\begin{proof}[Proof of Proposition~\ref{prop:reachcritical}]
Let $e$ be a unit eigenvector of $H$ with eigenvalue $\alpha$ of the sign of ${t}$. Since
$\phi(x_0+se)=c+\tfrac12\alpha s^{2}+o(s^{2})$, for small $|{t}|$ there is $s_{t}>0$ with $\phi(x_{t})=c+{t}$,
$x_{t}:=x_0+s_{t} e$, and $s_{t}=\sqrt{2{t}/\alpha}\,(1+o(1))$. Next, $\nabla\phi(x_{t})=s_{t}\alpha e+o(s_{t})$. Thus,
$|\nabla\phi(x_{t})|=\sqrt{2|{t}\alpha|}\,(1+o(1))$, and the unit normal tends to $\pm e$.

The compression of $H$ to $e^{\perp}$ has as eigenvalues the remaining eigenvalues of $H$, and there are some, as
$d\ge2$. Thus, its norm is at least $\alpha_{\min}$. By continuity of $D^{2}\phi$ and of the normal,
$\|P_nD^{2}\phi(x_{t})P_n\|\ge\alpha_{\min}(1-o(1))$, where $P_n$ projects onto the tangent space of $\Sigma_{t}$ at
$x_{t}$. With the formula for $\II$ in Proposition~\ref{prop:reachcritical}, this gives
$\|\II(x_{t})\|\ge\alpha_{\min}/\sqrt{3|{t}\alpha|}$ for small $|{t}|$. Finally,
$\rch(\Sigma_{t})\le1/\kappa_{\max}\le1/\|\II(x_{t})\|$ by the formula for the reach~\cite[Thm.~3.4]{Aamari2019}.
\end{proof}

\bibliographystyle{siamplain}
\bibliography{refs}
\end{document}